\documentclass[]{article}

\usepackage[margin=0.5in]{geometry}
\usepackage{graphicx}

\usepackage{color}
\usepackage{amsmath,amsfonts,amssymb,amsthm,bm,calrsfs,stmaryrd}
\usepackage{graphicx,natbib} 
\usepackage{subfigure} 
\usepackage{algorithm}
\usepackage{algorithmic}

\newtheorem{theorem}{Theorem}
\newtheorem{lemma}{Lemma}

\title{Approximate Joint Matrix Triangularization}
\author{Nicolo Colombo\thanks{LCSB, University of Luxembourg} \and Nikos Vlassis\thanks{Adobe Research, San Jose, CA}}
\begin{document}
\maketitle
\begin{abstract} 
We consider the problem of approximate joint triangularization of a set of noisy jointly diagonalizable real matrices.
Approximate joint triangularizers are commonly used in the estimation of the joint eigenstructure of a set of matrices, with applications in signal processing, linear algebra, and tensor decomposition.
By assuming the input matrices to be perturbations of noise-free, simultaneously diagonalizable ground-truth matrices, the approximate joint triangularizers are expected to be perturbations of the exact joint triangularizers of the ground-truth matrices.
We provide \emph{a priori} and \emph{a posteriori} perturbation bounds on the `distance' between an approximate joint triangularizer and its exact counterpart.  
The \emph{a priori} bounds are theoretical inequalities that involve functions of the ground-truth matrices and noise matrices, whereas the
\emph{a posteriori} bounds are given in terms of observable quantities that can be computed from the input matrices.
From a practical perspective, the problem of finding the best approximate joint triangularizer of a set of noisy matrices amounts to solving a nonconvex optimization problem.
We show that, under a condition on the noise level of the input matrices, it is possible to find a good initial triangularizer such that the solution obtained by any local descent-type algorithm has certain global guarantees.
Finally, we discuss the application of approximate joint matrix triangularization to canonical tensor decomposition and we derive novel estimation error bounds.
\end{abstract} 

\section{Introduction}
We address an estimation problem that appears frequently in engineering and statistics, whereby we observe noise-perturbed versions of a set of jointly decomposable matrices $M_n$, and the goal is to recover (within a bounded approximation) some aspects of the underlying decomposition. 
An instance of this problem is {\em approximate joint diagonalization}:
\begin{equation} 
\label{jointschur}
\hat M_n = M_n + \sigma W_n,  \qquad M_n = V^{} {\rm diag}([\Lambda_{n1}, \dots, \Lambda_{nd}]) V^{-1}, \qquad n=1,\ldots,N ,
\end{equation}
where $\hat M_n$ are the $d \times d$ observed matrices, and the rest of the model primitives are unobserved: $\sigma > 0$ is a scalar,
$W_n$ are arbitrary noise matrices with Frobenius norm $\| W_n \| \leq 1$, and  
the matrices $V, \Lambda$ define the joint eigenstructure of the ground-truth matrices $M_n$.
The optimization problem involves estimating from the observed matrices $\hat M_n$ the eigenvalues $\Lambda$ and/or the common factors $V$. 
Joint matrix diagonalization appears in many notable applications, such as independent component analysis \citep{cardoso1996jacobi}, latent variable model estimation \citep{balle2011spectral,AnimaJMLR14}, and tensor decomposition \citep{Lat06,Kuleshov15}.

Under mild conditions, the ground-truth matrices $M_n$ in \eqref{jointschur} can be \emph{jointly triangularized}, which is known as the (real) joint or simultaneous Schur decomposition \citep{horn2012matrix}. Namely, there exists an orthogonal matrix~$U_\circ$ that 
simultaneously renders all matrices $U_\circ^T M_n U_\circ$ upper triangular:
\begin{equation} 
\label{schur decomposition}
{\rm low}(U_\circ^T M_n U_\circ ) =  0  \qquad {\rm for \ all } \quad  n=1, \dots N,
\end{equation}
where low$(A)$ is the strictly lower-diagonal part of $A$ defined by 
$[\mbox{low}(A)]_{ij} = A_{ij}$ if $i>j$ and 0 otherwise.
On the other hand, when $\sigma>0$ the noisy matrices $\hat M_n$  in \eqref{jointschur} cannot be jointly triangularized exactly. The problem of {\em approximate joint triangularization} can be defined as the following optimization problem over the manifold of orthogonal matrices ${\mathbb O}(d)$:
\begin{equation}
\label{joint triangularization}
\min_{U \in {\mathbb O}(d)} {\cal L}(U) = \sum_{n=1}^N  \ \| {\rm low}(U^\top \hat M_n U) \|^2  \, .
\end{equation}
In words, we are seeking an orthogonal matrix $U$ such that the matrices $\hat T_n = U^T \hat M_n U$  are approximately upper triangular. 
This is a nonconvex problem that is expected to be hard to solve to global optimality in general. When $\sigma > 0$, the global minimum of ${\cal L}(U)$ will not be zero in general, and for any feasible $U \in {\mathbb O}(d)$ some of the entries below the main diagonal of each $\hat T_n$ may be nonzero.
The estimands of interest here could be the joint triangularizer $U$ and/or the approximate joint eigenvalues on the diagonals of $\hat T_n$.

Applications of (approximate) joint matrix triangularization range from algebraic geometry \citep{corless1997reordered}, to signal processing \citep{haardt_simultaneous_1998}, to tensor decomposition \citep{sardouie_canonical_2013,colombo2016icml}.
When the ground-truth matrices $M_n$ are symmetric, the models \eqref{jointschur} and \eqref{schur decomposition} are equivalent and $V,U_\circ$ are both orthogonal. However, when the matrices $M_n$ are non-symmetric, the matrix $V$ in \eqref{jointschur} is a general nonsingular matrix, while the matrix $U_\circ$ in \eqref{schur decomposition} is still orthogonal.
Since the optimization in \eqref{joint triangularization} is over a `nice' manifold, approximate joint triangularization is expected to be an easier problem than approximate joint diagonalization, the latter involving optimization over the manifold of invertible matrices \citep{afsari2008sensitivity}. 
Two types of methods have been proposed for optimizing~\eqref{joint triangularization}, Jacobi-like methods \citep{haardt_simultaneous_1998}, and Newton-like methods that optimize directly on the matrix manifold ${\mathbb O}(d)$ \citep{afsari2004some,colombo2016icml}.
Both methods are of iterative nature and their success depends on a good initialization.

\subsection{Contributions}
We are interested in theoretical guarantees for solutions $U$ computed by {\em arbitrary} algorithms that optimize \eqref{joint triangularization}.
Note that the objective function \eqref{joint triangularization} is continuous in the parameter $\sigma$.
This implies that, for $\sigma$ small enough, the approximate joint triangularizers of $\hat M_n$ can be expected to be perturbations of the exact triangularizers of $M_n$.
To formalize this, we express each feasible matrix $U$ in \eqref{joint triangularization} as a perturbation of an exact triangularizer $U_\circ$ of the ground-truth matrices $M_n$ in \eqref{jointschur}, that is
\begin{equation}
\label{expansion Ustar}
U  = U_\circ e^{\alpha X}, \quad {\rm where } \quad X = -X^\top ,  \quad \| X  \| = 1, \quad \alpha > 0 , 
\end{equation}
where $X$ is a skew-symmetric matrix and $e$ denotes matrix exponential.
Such an expansion holds for any pair $U,U_\circ$ of orthogonal matrices (see for example \citet{absil_optimization_2009}).
The scalar $\alpha$ in \eqref{expansion Ustar} can be interpreted as the `distance' between $U$ and $U_\circ$.

\paragraph{Perturbation bounds.}
We provide two different types of bounds on the parameter $\alpha$:
 \emph{A priori} bounds that are based on ground-truth quantities (such as the ground-truth matrices, the sample size, and in some applications also the assumed probability distribution generating the data), and {\em a posteriori} bounds that involve solely observable quantities (such as the observed matrices and the current solution).
While the former bounds are attractive theoretically as they can capture general perturbation effects on the matrix decomposition factors, the latter bounds can have more practical use, such as for instance in nonconvex optimization \citep{pang1987posteriori} and the design of optimized algorithms \citep{prudhomme2003practical}.

{\emph{A priori} analysis:} 
In Theorem \ref{theorem alpha} and Theorem \ref{theorem X} we provide 
two bounds that together offer a complete first-order characterization of the approximate triangularizers in terms of ground-truth quantities.
The corresponding inequalities depend on the noise level, the condition number of the joint eigenvectors matrix, a joint eigengap parameter, the number of ground-truth matrices, and their norm. Theorem \ref{theorem X} is the extension of the result derived by \citet{cardoso1994perturbation} for symmetric matrices. 

{\emph{A posteriori} analysis:} 
In Theorem \ref{theorem alpha star} we provide an error bound on the perturbation parameter $\alpha$, which is based on observable quantities that can be computed from the input matrices $\hat M_n$.
In particular, the bound involves the value of ${\cal L}(U)$ evaluated at each candidate $U$, and various functions of the observed matrices $\hat M_n$ and their approximate joint eigenvalues.
The only non-observable quantity appearing in the bound is the noise parameter $\sigma$ in \eqref{jointschur}, which, for instance in the case of joint matrix decomposition problems arising from empirical moment matching (see, e.g., \citet{AnimaJMLR14}), can be bounded by a function of the sample size.
The bound in Theorem \ref{theorem alpha star} is global, in the sense that it does not depend on the initialization, and can be used to characterize the output of any algorithm that optimizes \eqref{joint triangularization}.

\paragraph{Global guarantees for locally convergent algorithms.}
Beyond the purely theoretical analysis of approximate joint matrix triangularization, we also address the practical problem of computing an approximate joint triangularizer in \eqref{joint triangularization}.
Due to the nonconvexity of \eqref{joint triangularization}, locally convergent algorithms are guaranteed to converge to a given local optimum if the algorithm is started in the corresponding basin of attraction.
The continuity in the parameter $\sigma$  of the objective function ${\cal{L}}(U)$ in \eqref{joint triangularization} can be used to show that, under certain conditions, a finite set of local minimizers of \eqref{joint triangularization} enjoy global success guarantees in terms of their distance to the ground-truth matrices.  
In Theorem \eqref{theorem convergence} we provide a condition under which it is always possible to initialize a locally convergent algorithm in the basin of attraction of such a provably good minimizer.

\subsection{Related work}

The problem addressed here has two main antecedents: The work of  \citet{konstantinov1994nonlocal} on the perturbation of the Schur decomposition of a single matrix, and the work of \citet{cardoso1994perturbation} on the perturbation of joint diagonalizers. Our analysis can be viewed as an extension of the analysis of \citet{konstantinov1994nonlocal} to the multiple matrices case, and an extension of the analysis of \citet{cardoso1994perturbation} to joint matrix triangularization. We note that joint matrix triangularization is equivalent to joint spectral decomposition when the commuting matrices are symmetric. The proof of Theorem \ref{theorem X} exploits the same idea of \citet{cardoso1994perturbation}, but with a few key technical differences that pertain to non-symmetric / non-orthogonal matrices. We are not aware of other works dealing with the perturbation of joint matrix triangularizers. Moreover, to the best of our knowledge, our bound in Theorem \ref{theorem alpha star} is the first {\em a posteriori} error bound for joint matrix decomposition problems. 

From an algorithmic point of view, various approaches to approximate joint matrix triangularization have been proposed in the literature.
The simplest one is a matrix-pencil technique (see for example \citet{corless1997reordered}) where a linear combination of the input matrices is decomposed using established methods for the Schur decomposition of a single matrix.
The solution obtained in that case is, however, not optimal and depends on the particular matrix pencil.
A more standard way to formulate an approximate joint decomposition problem is to introduce a nonconvex objective function, as in \eqref{joint triangularization}, whose variables are the target shared matrix components \citep{cardoso1996jacobi,haardt_simultaneous_1998,abed1998least,fu2006balanced,Kuleshov15}.
The nonconvex optimization problem is then solved via iterative methods that typically belong to two classes, Jacobi-like methods \citep{cardoso1996jacobi,Kuleshov15}, and matrix manifold optimization methods \citep{afsari2004some,colombo2016icml}.
Jacobi-like algorithms rely on the decomposition of the variables into single-parameter matrices (such as Givens rotations), whereas
in a matrix manifold approach the objective \eqref{joint triangularization} is optimized directly on the matrix manifold.
As demonstrated recently \citep{colombo2016icml}, a Gauss-Newton method that optimizes \eqref{joint triangularization} directly on the matrix manifold ${\mathbb O}(d)$ can outperform the Jacobi-like method in terms of runtime by, roughly, one order of magnitude, for a statistically equivalent quality of the computed solutions.
Finally, the problem of obtaining global guarantees for joint matrix decomposition algorithms has been considered by \citet{Kuleshov15}, but only for the case of matrix joint diagonalization.
To the best of our knowledge, our work is the first that provides global solution guarantees for the joint matrix triangularization problem,
corroborating the strong empirical results that have been reported in the literature \citep{haardt_simultaneous_1998,abed1998least}.

\subsection{Conventions}
All matrices, vectors and numbers are real. Let $A$ be a $d \times d$  matrix, then $A^T$ is the transpose of $A$, $A^{-1}$ is the inverse of $A$ and $A^{-T}$ is the inverse of the transpose of $A$.
$A_{ij}$ (or $[A]_{ij}$) is the $(i,j)$ entry of $A$. 
The $i$th singular value of $A$ is denoted by  $\sigma_{i}(A)$  and $\kappa(A) = \frac{\sigma_{max}(A)}{\sigma_{min}(A)}$ is the condition number of $A$.
The matrix commutator $[A,B]$ is defined by $[A,B]=AB-BA$ and $\| A \| $ is the Frobenius norm defined by $\| A \|^2  = {\rm Tr}(A^T A) = \sum_{i,j} A_{ij}^2$.
The Kronecker product is denoted by $\otimes$.
Depending on the context, we use $1$ to denote a vector of ones or the identity matrix. 
${\mathbb O}(d)$ is the manifold of orthogonal matrices $U$ defined by $U^T U = 1$.
$T_{{\mathbb O}(d)}$ is the tangent space of ${\mathbb O}(d)$, \emph{i.e.} the set of skew-symmetric matrices satisfying $A = - A^T$.
${\rm vec}(A)$ is the column wise vectorization of $A$. 
${\rm low}(A)$ and ${\rm up}(A)$ are the strictly lower-diagonal and strictly upper-diagonal part of $A$ defined by
\begin{equation}
[ {\rm low}\left(A \right)]_{ij} = \left\{ 
\begin{array}{ll}  A_{ij} & {\rm if} \ i>j \\ 
0 & {\rm if} \ i \leq j  
\end{array} 
\right. 
\end{equation}
\begin{equation}
 [{\rm up}\left(A \right) ]_{ij} = \left\{ 
\begin{array}{ll}  A_{ij} & {\rm if} \ i<j \\
 0 & {\rm if} \ i \geq j  
 \end{array} 
 \right. 
\end{equation}
${\rm Low} \in \{0,1\}^{n^2\times n^2}$ and ${\rm Up}\in \{0,1\}^{n^2\times n^2}$ are linear operators defined by ${\rm vec}({\rm low}(A)) = {\rm Low} \ {\rm vec}(A)$ and ${\rm vec}({\rm up}(A)) = {\rm Up} \ {\rm vec}(A)$ respectively.
$P_{\rm Low} \in \{0,1\}^{\frac{n (n-1)}{2} \times n^2}$ is the projector to the sub-space of (vectorized) strictly lower-diagonal matrices defined by $P_{\rm Low} P_{\rm Low}^T = 1$ and  $P_{\rm Low}^T P_{\rm Low} = {\rm Low}$.
For example, letting $d=4$, one has
\begin{eqnarray}
{\rm Low} &=& {\rm diag}\left([0, 1, 1, 1,   0, 0, 1, 1,   0, 0, 0, 1,   0, 0, 0, 0]]\right) = {\rm diag}\left( 1^T P_{\rm low}^T \right) \\[5pt]
P_{\rm low} &=& \left( 
\begin{array}{cccc cccc cccc cccc}
0&1&0&0&  0&0&0&0&   0&0&0&0&   0&0&0&0\\ 
0&0&1&0&  0&0&0&0&   0&0&0&0&   0&0&0&0\\ 
0&0&0&1&  0&0&0&0&   0&0&0&0&   0&0&0&0\\ 
0&0&0&0&  0&0&1&0&   0&0&0&0&   0&0&0&0\\ 
0&0&0&0&  0&0&0&1&   0&0&0&0&   0&0&0&0\\ 
0&0&0&0&  0&0&0&0&   0&0&0&1&   0&0&0&0\\ 
\end{array}\right)
\end{eqnarray}
and similarly for ${\rm Up}$ and $P_{\rm up}$.

\section{Exact joint triangularizers}
Consider the set of simultaneously diagonalizable matrices ${\cal M}_{\circ} = \{ \hat M_n |_{\sigma = 0} \}_{n=1}^N$, with $\hat M_n$ defined in \eqref{jointschur}. 
A joint triangularizer of ${\cal M}_{\circ}$ is an orthogonal matrix $U_{\circ}$ such that 
\begin{equation}
{\rm low}( U_{\circ}^T M_n U_{\circ} ) = 0 \qquad {\rm for \ all } \quad n=1, \dots N  
\end{equation}
The condition under which ${\cal M}_{\circ}$ admits a finite number of joint triangularizers is established by the following lemma.
\begin{lemma}
\label{lemma non-degeneracy condition}
Let  ${\cal M}_{\circ} = \{ \hat M_n |_{\sigma = 0} \}_{n=1}^N$, with $\hat M_n$ defined in \eqref{jointschur}.
Then if 
\begin{equation}
\label{non-degeneracy condition lemma}
{\rm for \ every } \qquad  i \neq  i'  \qquad {\rm there \  exists } \qquad  n \in \{1, \dots, N \} \qquad {\rm s.t.} \quad \Lambda_{ni} \neq \Lambda_{ni'} 
\end{equation}
${\cal M}_{\circ}$ admits $2^d d!$ exact joint triangularizers.
\end{lemma}

\section{\emph{A priori} perturbation analysis}
Consider the approximate joint triangularization problem defined in \eqref{schur decomposition} and the expansion \eqref{expansion Ustar}. Two theoretical bounds are provided in this section. 
The first one is an inequality for the parameter $\alpha$.
The second one is an expression for the skew-symmetric matrix $X = -X^T$ that appears in \eqref{expansion Ustar}.
The explicit form of $X$ is given in terms of the ground-truth matrices $M_n$ and the noise matrices $\sigma W_n$.
Both bounds are valid up to second order terms in the perturbation parameters $\alpha $ and $\sigma$, \emph{i.e.} they hold up to $O((\alpha+\sigma)^2)$ terms.

\begin{theorem} 
\label{theorem alpha}
Let  ${\cal M}_{\sigma} =\{ \hat M_n \}_{n=1}^N $  and ${\cal M}_{\circ} = \{ \hat M_n |_{\sigma = 0} \}_{n=1}^N$ with $\hat M_n$ defined in \eqref{jointschur}.
Assume ${\cal M}_{\circ}$ is such that \eqref{non-degeneracy condition lemma} is satisfied. 
Then there exists $U_\circ$, which is an exact joint triangularizer of ${\cal M}_{\circ}$, such that an approximate joint triangularizer of ${\cal M}_{\sigma}$ can be written as 
\begin{equation}
\label{approximate triangularizer theorem}
U = U_{\circ} e^{\alpha X}  \qquad X = -X^T \qquad \|X \|  = 1  
\end{equation}
with $\alpha>0$ obeying 
\begin{equation}
\label{bound alpha theorem}
\alpha \leq  2 \sqrt 2 \sigma \| \tilde T^{-1}  \|_2   \sqrt{ \sum_{n=1}^N \| M_n\|^2}\sqrt{ \sum_{n=1}^N \| W_n\|^2}   + O((\alpha + \sigma)^2)
\end{equation}
where $M_n$ and $W_n$ are defined in \eqref{jointschur}, $\tilde T = \sum_{n=1}^N \tilde t_n^T \tilde t_n$ with $\tilde t_n = P_{\rm low} (1 \otimes  U_{\circ}^T M^T_n U_{\circ}  -  U_{\circ}^T M_n U_{\circ} \otimes ) P_{\rm low}^T$, and $\| \tilde T^{-1}\|_2$ is the spectral norm of the inverse of $\tilde T$.
\end{theorem}

It is possible to find a more explicit upper bound of \eqref{bound alpha theorem}, given  in terms of the ground  matrices and $\sigma$.
This result is provided by the following lemma
 
\begin{lemma}
\label{lemma bound non spectral}
Let $\alpha$ be defined as in Theorem \ref{theorem alpha}, then
\begin{equation}
\alpha \leq   \frac{ 2 \sigma \sqrt{d(d-1)} \kappa(V)^4}{\gamma} \sqrt{ \sum_{n=1}^N  \| M_n\|^2 } \sqrt{ \sum_{n=1}^N \| W_n\|^2}
 + O((\alpha+\sigma)^2)   \qquad \gamma = \min_{i < i'} \sum_{n=1}^N( \Lambda_{ni} -\Lambda_{ni'} )^2  
\end{equation}
where $V$, $M_n$, $W_n$ and $\Lambda$ are defined in \eqref{jointschur}.
\end{lemma}

\begin{theorem} 
\label{theorem X}
Let $U  =  U_{\circ} e^{\alpha X}$ be the approximate joint triangularizer defined in Theorem \ref{theorem alpha}. An approximate expression for the matrix $\alpha X$ is given by
\begin{eqnarray}
\alpha X = E -E^T \qquad E = {\rm mat}(P_{\rm low}^T x)  
\qquad 
\label{vec x first order}
 x  = - \sigma \left( \sum_{n=1}^N \tilde t_n^T \tilde t_n \right)^{-1}  \sum_{n=1}^N \tilde t^T_n P_{\rm low} {\rm vec}(U_{\circ}^T W_n U_{\circ}) + O((\alpha+\sigma)^2) 
\end{eqnarray}
where $ \tilde t_n = P_{\rm low} (1 \otimes U_{\circ}^T M^T_n U_{\circ} - U_{\circ}^T M_n U_{\circ} \otimes 1)P_{\rm low}^T $, with $M_n$ and $W_n$ defined in \eqref{jointschur}.   
\end{theorem}

\paragraph{Remarks on the theorems:}
The proof of these bounds is based on a first-order characterization of the approximate joint triangularizer $U$, which is defined as a stationary point of \eqref{joint triangularization}.
The inequalities on the parameter $\alpha$ come from the analysis  of the associated stationarity equation $\nabla {\cal L} = 0$, via a first order expansion around $U_\circ$, an exact joint triangularizer of ${\cal M}_{\sigma= 0}$.

\section{\emph{A posteriori} perturbation analysis}
\label{section a posteriori}
The result of this section is an {\em a posteriori} bound on the magnitude of the approximation error:

  \begin{theorem}
  \label{theorem alpha star}
Let  ${\cal M}_{\sigma} =\{ \hat M_n \}_{n=1}^N $  and ${\cal M}_{\circ} = \{ \hat M_n |_{\sigma = 0} \}_{n=1}^N$ with $\hat M_n$ defined in \eqref{jointschur}.
Assume that ${\cal M}_{\circ}$ satisfies \eqref{non-degeneracy condition lemma} and the noise matrices $W_n$ defined in \eqref{jointschur} obey $\| W_n\| \leq 1$.
Let $U$ be a feasible solution of the optimization problem \eqref{joint triangularization}.
Then there exists $U_\circ$, which is an exact joint triangularizer of ${\cal M}_{\circ}$, such that $U$ can be written as
\begin{equation}
\label{approximate triangularizer}
U = U_\circ e^{\alpha X} , \qquad X = -X^\top, \qquad \|X \|  = 1, \qquad \alpha >0  ,
\end{equation}
with $\alpha$ obeying
\begin{equation}
\label{alpha bound}
\alpha \leq  \sqrt 2 \| \beta \| \| \hat T_{\beta}^{-1} \|_2 (\sqrt{{\cal L}(U)} + \sigma \sqrt{N} )+ O((\sigma+\alpha)^2)
\end{equation} 
where $\beta = [\beta_1, \dots, \beta_N] \in {\mathbf R}^{N}$, $\hat T_{\beta} =  \sum_{n=1}^N \beta_n P_{\rm low}(1 \otimes U^T \hat M_n^T  U - U^T \hat M_n^T  U  \otimes 1) P_{\rm low}^T$, $\| T^{-1}_{\beta} \|_2 $ is the spectral norm of $T^{-1}_{\beta}$ and ${\cal L}(U)$ is defined in \eqref{joint triangularization}.
\end{theorem}

\paragraph{Remarks on the theorem:}
Assuming an a priori knowledge of $\sigma$, the inequality depends only on quantities that can be computed from the observed matrices $\hat M_n$.
The technique we have used to obtain the {\em a posteriori} bound follows an idea of \citet{konstantinov1994nonlocal} and is based on the perturbation equation
\begin{align}
\label{intro konstantinov}
U^\top \bigg( \sum_{n=1}^N \beta_n (M_n + \sigma W_n )\bigg)  U 
&= \sum_{n=1}^N \beta_n ( T_n + \varepsilon_n) , \qquad
{\rm low}(T_n) = 0 ,  \qquad &
\varepsilon_n = {\rm low} \bigg( U^\top\Big( \sum_{n=1}^N \beta_n \hat M_n \Big) U \bigg)
\end{align} 
where $\beta = [\beta_1, \dots, \beta_N]$.
The difference from the single matrix case studied by \citet{konstantinov1994nonlocal} is that the lower-diagonal terms $\varepsilon_n$ may be nonzero because an exact joint triangularizer may not exist.

\section{Global guarantees for locally convergent algorithms}
The existence of at least one approximate joint triangularizer of ${\cal M}_{\sigma} = \{ \hat M_n \}_{n=1}^N$ that is close to an exact triangularizer of ${\cal M}_{\circ} = \{ \hat M_n |_{\sigma=0}\}_{n=1}^N$ is guaranteed by the continuity of \eqref{joint triangularization} in the noise parameter $\sigma$.
The distance between such an approximate joint triangularizer, $U$, and the exact triangularizer, $U_{\circ}$, is bounded by the Theorem \ref{theorem alpha}.
If it is possible to compute a good initialization, a locally convergent algorithm is expected to converge to such $U$.
The following theorem provides a way to compute such a good initialization,  under certain conditions on the noise parameter $\sigma$. 

\begin{theorem}
\label{theorem convergence}
Let  ${\cal M}_{\sigma} =\{ \hat M_n \}_{n=1}^N $  and ${\cal M}_{\circ} = \{ \hat M_n |_{\sigma = 0} \}_{n=1}^N$ with $\hat M_n$ defined in \eqref{jointschur}.
Assume that ${\cal M}_{\circ}$ satisfies \eqref{non-degeneracy condition lemma} and the noise matrices $W_n$ defined in \eqref{jointschur} obey $\| W_n\| \leq 1$.
Let $\beta = [\beta_1, \dots, \beta_N] \in {\mathbf R}^N$  such that
\begin{equation}
\min_{i<i'} |{\rm Re}(\lambda_i(\hat M_\beta)-\lambda_{i'}(\hat M_\beta))| >0 \qquad \hat M_\beta=\sum_{n=1}^N \beta_n \hat M_n ,
\end{equation}
then a descent algorithm initialized with an orthogonal matrix  $U_{init}$ such that ${\rm low}(U_{init}^T \hat M_\beta U_{init}) = 0 $ (obtained via the Schur decomposition of $\hat M_\beta$) converges to an approximate joint triangularizer defined by Theorem \ref{theorem alpha} if the noise parameter $\sigma$ obeys
\begin{equation}
\sigma \leq \frac{2 \varepsilon }{ \sqrt{2N} \| \hat T^{-1}_\beta\|_2 A_{\alpha}    +  A_{\sigma} } + O(\sigma^2)
\end{equation}
where  
\begin{equation}
\varepsilon = \frac{\gamma}{2 \kappa(V)^4}  \qquad \gamma =   \min_{i<i'} \sum_{n=1}^N(\Lambda_{ni}-\Lambda_{ni'})^2 
\end{equation}
\begin{equation}
 \hat T_\beta  = P_{\rm low} (1 \otimes U^T_{init} \hat M^T_\beta U_{init} - U^T_{init} \hat M_\beta U_{init} \otimes 1) P_{\rm low} \qquad A_\alpha = 32  \sum_{n=1}^N \| M_n\|^2 \qquad A_\sigma = 16  \sqrt{N} \sqrt{\sum_{n=1}^N \| M_n \|^2}
\end{equation}
with $M_n$, $V$ and $\Lambda$ defined in \eqref{jointschur}.
\end{theorem}   
\paragraph{Remarks on the theorem:}
The proof of the theorem consists of two steps: 

(i) We first characterize the convex region containing an exact joint triangularizer $U_{\circ}$, in terms of $\alpha_{max}$,  the distance from $U_{\circ}$. This is obtained by requiring that the Hessian of \eqref{joint triangularization} computed at $U = U_{\circ} e^{\alpha X} $ is positive definite for all $X$ (with $\|X \|=1 $) if $\alpha \leq \alpha_{\max}$. 

(ii) Then we find a condition on the noise parameter for which the orthogonal matrix $U_{init}$, which is used to initialize the algorithm, belongs to the convex region characterized in the previous step. Letting $U_{init}=U_{\circ} e^{\alpha_{init}X_{init}}$, this is equivalent to $\alpha_{init} \leq \alpha_{max}$. 

Global success guarantees for the solution $U$ computed by a local hill-climbing algorithm can be obtained by combining Theorem~\ref{theorem convergence} and Theorem \ref{theorem alpha}.

\section{Applications to tensor decomposition}
\subsection{Observable matrices}
Consider an order the $N \times N \times N$ tensor of the form
\begin{equation}
\label{tensor}
\hat {\mathbb T} = {\mathbb T} + \sigma {\mathbb E} \qquad {\mathbb T}_{nn' n''} =  \sum_{i=1}^d Z_{ni} Z_{n'i} Z_{n''i} \qquad n,n',n'' = 1,\dots N \  
\end{equation}
where $\sigma >0$ and ${\mathbb E}$ is an arbitrary noise term satisfying $\|{\mathbb E} \|\leq \varepsilon$,  with $\|{\mathbb E} \|^2 = \sum_{n n' n''}{\mathbb E}_{n n' n''}^2$.
We define the $d \times d$ `observable' matrices associated with the tensor $\hat {\mathbb T}$ as
\begin{eqnarray}
\label{observable matrices general R}
\hat M_{n}= \hat m_n  \ \hat m^{-1} \qquad n=1,\dots, N \qquad \hat m = \sum_{n=1}^N \hat m_n 
\end{eqnarray}
where, for general $d\leq N$, $\hat m_n$ are dimension-reduced tensor slices defined by
\begin{eqnarray}
\hat m_n =  U_d^T  \tilde m_n V_d \qquad   [\tilde m_n]_{n'n''}=\hat {\mathbb T}_{n n' n''} \qquad n,n',n''=1,\dots N   
\end{eqnarray}
 with $U_d$ and $V_d$ being $N\times d$ Stiefel matrices obtained by staking the first $d$ left and right singular vectors of $\sum_{n=1}^N \tilde m_n$.
The definition \eqref{observable matrices general R} makes sense only if $\hat m$ is invertible, \emph{i.e.} if the $d$th singular value of $\sum_{n=1}^N \tilde m_n$ is non-vanishing.
Assuming $d = N$ there is no need of introducing the dimension reduction matrices $U_d$ and $V_d$ and the observable matrices are then defined by 
\begin{equation}
\label{observable matrices}
\hat M_{n}= \hat m_n  \ \hat m^{-1} \qquad [\hat m_n]_{n' n''}=\hat {\mathbb T}_{n n' n''} \qquad n,n',n''=1,\dots, N \quad (d = N) \qquad \hat m = \sum_{n=1}^N \hat m_n  
\end{equation}
where $\hat m$ is assumed to be invertible. A more general definition of $\hat m$ would be $\hat m_{\theta} = \sum_{n=1}^N \theta_n \hat m_n$ where $\theta$ is an arbitrary $N$-dimensional vector. In what follows we consider the case $d=N$ and $\theta = 1$ but generalizations to $d\leq N $ and $\theta \neq 1$ are straightforward.
Observable matrices of the form \eqref{observable matrices general R} cannot be defined if $d>N$.
Given \eqref{tensor} and \eqref{observable matrices}, it is easy to prove the following lemma
\begin{lemma}
\label{lemma observable matrices}
If $Z$ is invertible and $[1^T Z]_i \neq 0$ for all $i=1, \dots, d$, the observable matrices $\hat M_n$ defined in \eqref{observable matrices} can be expanded as follows
\begin{eqnarray}
\label{observable matrices expansion}
\hat M_{n} = M_n + \sigma W_n + O(\sigma^2) \qquad M_n = Z {\rm diag}({\bf e}_n^T Z) \left( {\rm diag}(1^T Z)\right)^{-1} Z^{-1}  \qquad W_n = e_n m^{-1} + m_n m^{-1} e m^{-1}
\end{eqnarray}
where $ n= 1,\dots, N$, the vector ${\bf e}_n$ is the $n$th basis vector, and
\begin{equation}
 [e_n]_{n'n''}=E_{nn'n''} \qquad e=\sum_{n=1}^{N} e_n \qquad [m_n]_{n'n''}={\mathbb T}_{nn'n''} \qquad m = \sum_{n=1}^{N} m_n \, .
 \end{equation}
If ${\mathbb E}$ in \eqref{tensor} obeys $\| {\mathbb E} \| \leq \varepsilon $, then 
\begin{equation}
\| M_n \|  \leq   \frac{ d  \ \kappa(Z)^2 \max | Z |}{\min | 1^T Z |} \qquad 
\| W_n \|  \leq  \frac{\varepsilon \ \kappa(Z)^2 \  \sqrt{d}}{\| Z \|^2 \min | 1^T Z |} \left(1 +  \frac{ d  \ \kappa(Z)^2 \max | Z |}{\min | 1^T Z |} \right)  
\end{equation}
\end{lemma}

\subsection{Estimation of the tensor components $Z$} 
Lemma \ref{lemma observable matrices} implies that $Z$ can be obtained, up to normalization constants, from the estimated joint eigenvalues of the nearly jointly diagonalizable matrices \eqref{observable matrices expansion}.
Let $U$ be an approximate joint triangularizer of ${\cal M}_{\sigma} = \{ \hat M_n \}_{n=1}^N$ obeying the bound in Theorem \ref{theorem alpha}.
The corresponding estimation of $Z$ is given by   
\begin{equation}
\label{estimation Z}
\frac{Z^*_{ni}}{[1^T Z^*]_i}  =  [ U^T  \hat M_n  U  ]_{ii} \qquad n=1,\dots,N \qquad i=1,\dots,d 
\end{equation} 
where $[1^T Z^*]_i$ is an undetermined column-rescaling factor and we assume $ N = d$.
Under the conditions that $Z$ is invertible and $[1^T Z]_i\neq 0$ for all $i=1,\dots,d$, the difference between the estimated tensor components \eqref{estimation Z} and the ground-truth tensor components $Z$ is bounded by the following theorem. 
\begin{theorem}
\label{main theorem}
Let $\hat {\mathbb T}$ be the tensor defined in \eqref{tensor} and assume $N=d$,  
$Z$ is  invertible, and $[1^T Z]_i \neq 0$ for all $i=1,\dots d$.
Let $U$ be an approximate joint triangularizer of ${\cal M}_{\sigma} = \{ \hat M_n \}_{n=1}^N$, with $\hat M_n$ defined in \eqref{observable matrices expansion}, and $\frac{Z^*_{ni}}{[1^T Z^*]_i} = [ U^T \hat M_n U]_{ii}$ for all $n=1, \dots N$ and $i=1, \dots d$.
Then, if $U$ obeys the bound in Theorem \ref{theorem alpha}, $Z_*$ is such that
\begin{eqnarray}
\label{main bound}
\left| \frac{Z_{ni}^*}{[1^TZ^*]_i} - \frac{Z_{ni}}{[1^TZ]_i} \right|  \leq 4\sigma  \frac{\sqrt{d(d-d)} \kappa(Z)^4}{\gamma} {\rm M}^2 {\rm W} + \sigma {\rm W} + O(\sigma^2)
\end{eqnarray}
where 
\begin{equation}
\gamma = \frac{1}{N}\min_{i\neq i'} \sum_{n=1}^{N} (Z_{n i}-Z_{n i'})^2  \qquad
{\rm M} \leq \frac{ N   \kappa(Z)^2   \max | Z |}{\min | 1^T Z |}\qquad 
{\rm W} \leq   \frac{\varepsilon    \sqrt{N}   \kappa(Z)^2 }{\| Z \|^2 \min | 1^T Z |} \left(1 +   \frac{ N \kappa(Z)^2 \max | Z |}{\min | 1^T Z |} \right) 
\end{equation}
\end{theorem}

\paragraph{Remark on the theorem:}
Theorem \ref{main theorem} provides a first order approximation of the estimation error and it is valid up to terms proportional to $\sigma^2$. 
The assumption on $[1^T Z]_i$ can be relaxed by defining $\hat m$ as $\hat m_{\theta} = \sum_{n=1}^N \theta_n \hat m_n $, where $\theta$ is any $N$-dimensional vector for which $[\theta^T Z]_i\neq 0$ for all for $i=1,\dots, d$. 
The normalization constants $[\theta^T Z^*]_i$ can then be obtained from $\hat m_{\theta}  = \sum_{n=1}^N \theta_n \hat m_n $ and the corresponding estimates $\frac{Z^*_{ni}}{[\theta^T Z]_i}$ by solving the following matrix equation 
\begin{eqnarray}
\hat m_{\theta}  =   Z^*\frac{1}{{\rm diag}(\theta^T Z^*)} {\rm diag}(\theta^T Z^*)^3 \left(Z^* \frac{1}{{\rm diag}(\theta^T Z^*)} \right)^T
\quad \qquad \left[ Z^* \frac{1}{{\rm diag}(\theta^T Z^*)} \right]_{ni} =[U^T \hat M_n U ]_{ii}
\end{eqnarray}
Finally, by using \eqref{estimation Z} and the \emph{a posteriori} error analysis  of Section \ref{section a posteriori} it is possible to obtain analogous bounds that depend only on the observable matrices \eqref{observable matrices}.

\section{Other lemmas and proofs}

\subsection{Proof of Lemma \ref{lemma non-degeneracy condition}}
Lemma \ref{lemma non-degeneracy condition} establishes a sufficient condition for the existence of $2^d d!$ exact joint triangularizers of ${\cal M}_{\circ} = \{ \hat M_n |_{\sigma=0}\}_{n=1}^N$, with $\hat M_n$ defined in \eqref{jointschur}.    
The proof consists of showing that, if \eqref{non-degeneracy condition lemma} holds (i) there exist $2^d d!$ exact joint triangularizers of t${\cal M}_{\circ}$ and (ii) it is impossible to find more than $2^d d!$ such orthogonal matrices. 
Lemma \ref{lemma linear combination} can be used to prove that, when \eqref{non-degeneracy condition lemma} is fulfilled, it is possible to define a linear combination of the matrices $M_n \in {\cal M}_{\circ}$ with distinct eigenvalues.
Let $M$ be such linear combination of the matrices $M_n$.
Since any real $d\times d$ matrix with distinct eigenvalues admits $2^d d!$ triangularizers, $M$ admits $2^d d!$ triangularizers.
Now, since $[M_n,M_{n'}] = 0$ one has 
\begin{equation}
\label{commutation linear combination}
[M, M_n] = 0  \qquad \forall  \ n=1,\dots,N
\end{equation}
implying that all $2^d d!$ triangularizers of $M$ exactly triangularize all $M_n \in {\cal M}_{\circ}$.
This is due to the fact that commuting matrices are always joint triangularizable and implies that ${\cal M}_{\circ}$ has at least $2^d d!$ joint triangularizers.
But the commutation relation \eqref{commutation linear combination} also implies that any possible additional triangularizer of a matrix $M_n \in {\cal M}_{\circ}$ would exactly triangularize $M$.
This contradicts the fact that $M$ admits only $2^d d!$ exact triangularizers and proves the lemma.
$\square$
\subsection{Proof of Theorem \ref{theorem alpha} }
The stationary point of \eqref{joint triangularization} are defined by the equation $\nabla {\cal L} = 0$ where $\nabla {\cal L}$ is the gradient of ${\cal L}$ and ${\cal L}$ is defined in \eqref{joint triangularization}.
According to Lemma \ref{lemma stationarity equation}, if $U $ is a stationary point of \eqref{joint triangularization}, then
\begin{equation}
\label{stationarity equation}
\nabla {\cal L} = S - S^T   = 0 \qquad S = \sum_{n=1}^N \left[ U^T \hat M^T_n U , {\rm low}( U^T \hat M_n U)  \right]
\end{equation}
Now, let $U = U_{\circ} e^{\alpha X}$, where $U_{\circ}$ is an exact triangularizer of ${\cal M}_{\circ} = \{ \hat M_n|_{\sigma=0}\}_{n=1}^N$, $\hat M_n$ are defined in \eqref{jointschur}, $X = -X^T$ and one can assume $\|X \|  = 1$ and $\alpha>0$.
The expansion of $S$ in $\alpha$ and $\sigma$ reads 
\begin{eqnarray}
S &=& S|_{(\alpha=0,\sigma=0)} + \alpha \partial_\alpha S|_{\sigma = 0} + \sigma \partial_\sigma S|_{\alpha=0} + O((\alpha + \sigma)^2) \\
& = & \sum_{n=1}^N \left[ U_{\circ}^T M^T_n U_{\circ} , {\rm low}( [U_{\circ}^T M_n U_{\circ}, \alpha X ] )  \right] +  \sum_{n=1}^N \left[ U_{\circ}^T M^T_n U_{\circ} , {\rm low}( U_{\circ}^T \sigma W_n U_{\circ} )  \right] + O((\alpha + \sigma)^2)
\end{eqnarray}
where we have defined $\partial_\alpha f = \frac{\partial}{\partial \alpha} f |_{\alpha = 0}$ and $\partial_\sigma f = \frac{\partial}{\partial \sigma} f |_{\sigma = 0}$.
Note that, for all $n=1, \dots, N$,  $ [ U_{\circ}^T M^T_n U_{\circ} , {\rm low}( A)] $ is strictly lower-triangular  for any $A$
because $ {\rm up}(U_{\circ}^T M^T_n U_{\circ}) = 0$.
The latter follows from the fact that $U_{\circ}$ is an exact triangularizer of ${\cal M}_{\circ}$ and hence $ U_{\circ}^T M_n U_{\circ}$ is upper triangular, for all $n=1, \dots, N$.
Considering only the lower-diagonal part of the stationarity equation one obtains the necessary condition
\begin{equation}
\label{necessary condition}
0 = {\rm low}(S - S^T) = {\rm low}(\alpha \partial_\alpha S|_{\sigma = 0} + \sigma \partial_\sigma S|_{\alpha=0} ) + O((\alpha + \sigma)^2))
\end{equation} 
since the first order terms of $S^T$ are upper triangular.
The projected stationarity equation \eqref{necessary condition} reads
\begin{eqnarray}
\label{S alpha}
{\rm low}\left( \sum_{n=1}^N \left[ U_{\circ}^T M^T_n U_{\circ} , {\rm low}( [U_{\circ}^T M_n U_{\circ}, \alpha X ] )  \right] \right)  = - {\rm low}\left(\sum_{n=1}^N \left[ U_{\circ}^T M^T_n U_{\circ} , {\rm low}( U_{\circ}^T W_n U_{\circ} )  \right]  \right) + O((\alpha+\sigma)^2) 
\end{eqnarray}
Moreover, since ${\rm low}(U_{\circ}^T M_n U_{\circ}) = 0$ for all $n=1, \dots, N$ one has 
\begin{equation}
{\rm low}( [U_{\circ}^T M_n U_{\circ}, \alpha X ] ) = {\rm low}( [U_{\circ}^T M_n U_{\circ}, {\rm low}(\alpha X) ] )  
\end{equation}
This means that the linear operator defined by
\begin{equation}
{\cal T} {\rm low} (X)  = {\rm low} \left( \sum_{n=1}^N \left[ U_{\circ}^T M^T_n U_{\circ} , {\rm low}( [U_{\circ}^T M_n U_{\circ}, {\rm low} (X) ] )  \right] \right)
\end{equation}
maps the subspace of strictly lower dimensional matrices into itself.
This is a $\frac{d(d-1)}{2}$-dimensional subspace that has the same degrees of freedom as the set of $d \times d$ skew-symmetric matrices.
Each  $d \times d$ skew-symmetric matrix is mapped into this subspace by means of the projection $P_{\rm low} {\rm vec}(X)$. 
Conversely, letting $x$ be an element of this subspace, the corresponding  $d \times d$ skew-symmetric matrix $X$ is given by $X = {\rm mat}(P_{\rm low}^T x) - {\rm mat}(P_{\rm low}^T x)^T$.
Let $T$ be the linear operator defined by the vectorization of \eqref{S alpha}
\begin{equation}
T = \sum_{n=1}^N t_n^T t_n \qquad t_n = {\rm Low} (1 \otimes  U_{\circ}^T M^T_n U_{\circ}  -  U_{\circ}^T M_n U_{\circ} \otimes ){\rm Low}
\end{equation}
Its reduction to the subspace of strictly lower-diagonal matrices can be written as
\begin{equation}
\tilde T  = P_{\rm low}  T P_{\rm low}^T = \sum_{n=1}^N \tilde t_n^T \tilde t_n \qquad \tilde t_n = P_{\rm low} (1 \otimes  U_{\circ}^T M^T_n U_{\circ}  -  U_{\circ}^T M_n U_{\circ} \otimes )P_{\rm low}^T
\end{equation}
Then one has 
\begin{equation}
\label{linearization S alpha}
P_{\rm low}{\rm vec}({\cal T} {\rm low} (\alpha X)) = \tilde T P_{\rm low}{\rm vec}(\alpha X) 
\end{equation}
The $\frac{d(d-1)}{2} \times \frac{d(d-1)}{2}$ $\tilde T$ is positive definite if the non-degeneracy condition in \eqref{non-degeneracy condition lemma} is fulfilled (see Lemma \ref{lemma T}).
Under this assumption  
\begin{equation}
\label{solution X}
\alpha P_{\rm low} {\rm vec}(X)  = - \tilde T^{-1}  \ P_{\rm low} {\rm vec}\left(\sum_{n=1}^N \left[ U_{\circ}^T M^T_n U_{\circ} , {\rm low}( U_{\circ}^T W_n U_{\circ} )  \right]  \right) + O((\alpha+\sigma)^2) 
\end{equation}
Taking the norm of both sides one has
\begin{eqnarray}
\label{final eq alpha}
\alpha \leq 2 \sqrt{2}  \sigma  \| \tilde T^{-1} \|_2 \sqrt{ \sum_{n=1}^N  \| M_n\|^2 } \sqrt{ \sum_{n=1}^N  \| W_n\|^2 } 
\end{eqnarray} 
where we have used  $\|{\rm low}(X)\| = \frac{1}{\sqrt 2} \| X\|$, $\| X \|  = 1$ and 
\begin{eqnarray}
\left\| P_{\rm low} {\rm vec}\left(\sum_{n=1}^N \left[ U_{\circ}^T M^T_n U_{\circ} , {\rm low}( U_{\circ}^T W_n U_{\circ} )  \right]  \right) \right\|  \leq  2 \sigma \sqrt{ \sum_{n=1}^N \| M_n\|^2 }  \sqrt{ \sum_{n=1}^N \| W_n\|^2 } 
\end{eqnarray}
from $\|  \sum_{n=1}^N  t^T_n {\rm vec}( U_{\circ}^T \sigma W_n U_{\circ} )  \| \leq  \sqrt{ \sum_{n=1}^N  \|t^T_n \|^2 }\sqrt{ \sum_{n=1}^N  \| \sigma W_n\|^2} $, $\| t_n\|^2 \leq 4 \| M_n\|^2$.
$\square$

\subsection{Proof of Lemma \ref{lemma bound non spectral}}
Consider the inequality on the perturbation parameter $\alpha$ given in \eqref{final eq alpha}.
Lemma \ref{lemma T} states that the matrix $\tilde T$ is positive definite if the non-degeneracy condition \eqref{non-degeneracy condition lemma} is fulfilled and in this case \begin{equation}
\| \tilde T^{-1} \|_2 \leq  \sqrt{\frac{d(d-1)}{2} }  \frac{\kappa(V)^4}{\gamma}  \qquad \gamma = \min_{i<i'} \sum_{n=1}^N (\Lambda_{ni} - \Lambda_{ni'})^2
\end{equation}
This implies 
\begin{equation}
\alpha \leq     \frac{2 \sigma \sqrt{d(d-1)} \kappa(V)^4}{\gamma} \sqrt{ \sum_{n=1}^N  \| M_n\|^2 } \sqrt{ \sum_{n=1}^N  \| W_n\|^2 }
\end{equation}
$\square$ 
\subsection{Proof of Theorem \ref{theorem X} }
Theorem \ref{theorem X} follows from \eqref{solution X} where one can use 
\begin{eqnarray}
P_{\rm low} {\rm vec} \left(\sum_{n=1}^N \left[ U_{\circ}^T M^T_n U_{\circ} , {\rm low}( U_{\circ}^T \sigma  W_n U_{\circ} )  \right] \right) =  \sigma \sum_{n=1}^N \tilde t^T_n P_{\rm low}(U_{\circ}^T W_n U_{\circ}) 
\end{eqnarray} 
to obtain
\begin{equation}
P{\rm low} {\rm vec}(\alpha X)  = - \sigma  \left(\sum_{n=1}^N \tilde t^T_n \tilde t_n \right)^{-1}  \sum_{n=1}^N \tilde t^T_n P_{\rm low}(U_{\circ}^T W_n U_{\circ})
\end{equation}
with $\tilde t_n = P_{\rm low} (1 \otimes  U_{\circ}^T M^T_n U_{\circ}  -  U_{\circ}^T M_n U_{\circ} \otimes 1 ) P_{\rm low}^T$.
$\square$

\subsection{Proof of Theorem \ref{theorem alpha star}}
Let $\sum_{n=1}^N \beta_n \hat M_n$ be a general linear combination of the input matrices, where $\beta_n$, $n=1, \dots, N$ are arbitrary real numbers.
Let $U_\circ$ be an exact joint triangularizer of ${\cal M}_{\circ}$,  and $U$ be a feasible solution of the joint triangularization problem \eqref{joint triangularization}.
By construction $U$ is an orthogonal matrix and can be written as $U = U_\circ e^{\alpha X}$, with $X = -X^\top $, $\|X \| = 1$ and $\alpha >0$.
For any choice of $\beta$ one has
\begin{equation}
U^\top \left( \sum_{r=1}^N  \beta_n \hat M_n \right) U = \sum_{r=1}^N  \beta_n (\hat T_n + \varepsilon_n) \qquad {\rm low}(\hat T_n) = 0 \quad \varepsilon_n = {\rm low}(U^\top \hat M_n U)  
\end{equation}
By projecting onto the strictly lower-diagonal part and considering the expansion $U = U_\circ e^{\alpha X}$, we obtain 
\begin{eqnarray}
 \sum_{r=1}^N  \beta_n \varepsilon_n & = &  \sum_{r=1}^N \beta_n {\rm low}\left(e^{ - \alpha X} U_\circ^\top M_n U_{\circ}e^{\alpha X} +  e^{ - \alpha X} U_\circ^\top \sigma W_n U_\circ e^{\alpha X} \right)  \\
\label{equation konstantinov}
& = &\sum_{r=1}^N \beta_n  {\rm low}\left( [U_\circ^\top M_n U_\circ, \alpha X] +  U_\circ^\top \sigma W_n U_\circ \right) + O((\alpha + \sigma)^2)
\end{eqnarray}   
For any $X$, one has ${\rm low}([U_\circ^\top M_n U_\circ, X]) ={\rm low}([ U_\circ^\top M_n U_\circ,{\rm low}(X)])$ 
because $U_\circ^\top M_n U_\circ$ is upper triangular.  
The identity \eqref{equation konstantinov} can be rewritten as 
\begin{equation}
{\rm low}\left([U_\circ^\top \sum_{r=1}^N \beta_n  M_n U_\circ, {\rm low}(\alpha X)]\right)  = \sum_{r=1}^N \beta_n {\rm low}\left( \varepsilon_n  - U_\circ^\top \sigma W_n U_\circ \right)
\end{equation}
whose vectorization reads
\begin{equation}
\label{vectorized equation}
T_{\beta} {\rm vec}(\alpha X)  = {\rm vec}\left( {\rm low}\left(\sum_{r=1}^N \beta_n\varepsilon_n  -   \sigma W_{\beta}\right) \right) \qquad T_{\beta} = {\rm Low}(1 \otimes M_{\beta}^\top - M_{\beta} \otimes 1) {\rm Low}
\end{equation}
where  $M_{\beta} = \sum_{r=1}^N \beta_n U_\circ^\top M_n U_\circ $ and $W_{\beta} = \sum_{r=1}^N \beta_n U_\circ^\top W_n U_\circ$.
The reduction of $T_{\beta}$ to the subspace of strictly lower-diagonal matrices is
\begin{equation}
\tilde T_{\beta} = P_{\rm low} T_{\beta} P_{\rm low} ^T  = P_{\rm low} (1 \otimes M_{\beta}^\top - M_{\beta} \otimes 1) P_{\rm low}^T
\end{equation}
Lemma \ref{lemma T beta} can be used to show that $\tilde T_\beta$ is invertible if $M_{\beta}$ is invertible and $\lambda_{i}(M_\beta) \neq \lambda_{i'}(M_\beta)$ for all $i\neq i'$. 
Under this assumption one can write
\begin{equation}
\alpha  {\rm vec}(X)  =  \tilde T_{\beta}^{-1} {\rm vec}\left( {\rm low}\left(\sum_{r=1}^N \beta_n\varepsilon_n  -   \sigma W_{\beta}\right) \right) + O((\alpha + \sigma)^2) 
 \end{equation}
 and, by taking the norm in both sides, 
 \begin{eqnarray}
\alpha  &\leq&  \sqrt 2 \| \tilde T_{\beta}^{-1}\|_2  \| \beta \| \left( \sqrt{ \sum_{r=1}^N \| \varepsilon_n \|^2 }  +   \sigma  \sqrt{ \sum_{n=1}^N \| W_n\|^2 } \right) + O((\alpha + \sigma)^2) \\
& \leq &  \sqrt 2 \| \tilde T_{\beta}^{-1}\|_2 \left( \sqrt{ {\cal L}(U)} +   \sigma  \sqrt N\right) + O((\alpha + \sigma)^2) 
 \end{eqnarray}
 where we have used the assumption $\|\beta \|=1$ and $\| W_n \| \leq 1$.
Finally, one has 
\begin{eqnarray}
 \tilde T_{\beta}  &=& \sum_{r=1}^N \beta_n P_{\rm low}(1 \otimes U_{\circ}^T M_n^T  U_{\circ} - U_{\circ}^T M_n^T  U_{\circ}  \otimes 1) P_{\rm low}^T \\
& = &    \sum_{r=1}^N \beta_n P_{\rm low}(1 \otimes U_{\circ}^T \hat M_n^T  U_{\circ} - U_{\circ}^T \hat M_n^T  U_{\circ}  \otimes 1) P_{\rm low}^T  + O(\sigma)\\
& = &   \sum_{r=1}^N \beta_n P_{\rm low} (1 \otimes U^T \hat M_n^T  U - U^T \hat M_n^T  U  \otimes 1) P_{\rm low}  + O(\sigma+\alpha) \\
& = &  \hat T_{\beta} + O(\sigma+\alpha)
\end{eqnarray}
where we have defined $\hat T_\beta = \sum_{r=1}^N \beta_n P_{\rm low} (1 \otimes U^T \hat M_n^T  U - U^T \hat M_n^T  U  \otimes 1) P_{\rm low}$. 
It follows that $\| \tilde T^{-1}_{\beta} \|_2  = \|\hat T^{-1}_{\beta}\|_2 + O(\sigma+\alpha)$ and hence 
\begin{eqnarray}
\alpha  
& \leq &  \sqrt 2 \| \hat T_{\beta}^{-1}\|_2 \left( \sqrt{ {\cal L}(U)} +   \sigma  \sqrt N\right) + O((\alpha + \sigma)^2)
\end{eqnarray}
$\square$

\subsection{Proof of Theorem \ref{theorem convergence} }
The Hessian of \eqref{joint triangularization} at $U$ is positive definite if, for all $X$ such that $X = -X^T$,  $\langle X , \nabla^2 {\cal L} X \rangle >0$, where  
\begin{equation}
\langle X , \nabla^2 {\cal L}(U) X \rangle = \frac{d^2}{dt^2} {\cal L}(Ue^{t X})|_{t=0} 
\end{equation}
Lemma \ref{lemma hessian} shows that this is the case if 
\begin{equation}
\label{alpha max}
U = U_{\circ} e^{\alpha Y} \qquad Y = -Y^T \qquad \| Y \| = 1 \qquad  \alpha\leq \alpha_{max} \qquad \alpha_{max} = \frac{2 \varepsilon - \sigma A_{\sigma}}{A_\alpha} + O((\alpha+\sigma)^2)
\end{equation} 
\begin{equation}
\varepsilon = \frac{\gamma}{2 \kappa(V)^4}  \qquad \gamma = \min_{j<j'} \sum_{n=1}^N(\Lambda_{nj}-\Lambda_{nj'})^2
\qquad 
A_\alpha = 32 \sum_{n=1}^N \| M_n\|^2  \qquad 
A_\sigma = 16 \sqrt{N} \sqrt{\sum_{n=1}^N \| M_n \|^2}
\end{equation}
where we have assumed $\| W_n \| \leq 1$. 
The condition under which the Hessian of \eqref{joint triangularization} at $U_{\circ}$ is positive definite is $\alpha_{max}>0$.
If $U$ is a minimizer of ${\cal L}(U)$, this condition ensures that $U_{\circ}$ belongs to the convex region centered in $U$.
Now, assume that it is possible to find a vector $\beta = [\beta_1, \dots , \beta_N]$ such that $\| \beta\|=1$ and the operator  $T_{\beta}$ defined by
\begin{equation}
\label{initial T}
T_{\beta} = P_{\rm low} (1 \otimes U^T_{init} \hat M^T_\beta U_{init} - U^T_{init} \hat M_\beta U_{init} \otimes 1) P_{\rm low} \qquad U_{int} \in {\mathbb O}(d) \ {\rm s.t.} \   {\rm low}(U_{int}^T \hat M_\beta U_{int}) = 0 \qquad \hat M_\beta = \sum_{n=1}^N \beta_n \hat M_{n} 
\end{equation}
is invertible.
The orthogonal matrix $U_{int} $ is defined by the Schur decomposition of $\hat M_\beta$. According to Lemma \ref{lemma T beta}, $T_{\beta} $ is invertible if $\hat M_\beta$ is invertible and has real separated eigenvalues, \emph{i.e.} if $\lambda_i(\hat M_\beta)$ are real for all $i=1, \dots, d$ and $\min_{i<i'} |\lambda_i(\hat M_\beta) - \lambda_i(\hat M_\beta)| >0$. 
Finding such a $\hat M_\beta$ is possible if $\sigma$ is small enough. 
This is a consequence of Lemma \ref{lemma linear combination} and standard eigenvalues perturbation results. 
Otherwise, the separation of the eigenvalues of $\hat M_\beta$ can be checked numerically, since $\hat M_\beta$ is an observable quantity.
Now, let $M_\beta = \sum_{n=1}^N \beta_n M_n$, $W_\beta = \sum_{n=1}^N \beta_n W_n$  and $U_{\circ}\in {\mathbb O}(d) $ be such that ${\rm low}(U_{\circ}^T M_\beta U_{\circ}) = 0$.
By writing $U_{\circ} = U_{init} e^{\alpha Y}$ one has  
\begin{equation}
U_{\circ}^T M_\beta U_{\circ} = e^{-\alpha Y} U^T_{init} (\hat M_\beta - \sigma W_\beta) U_{init} e^{\alpha Y}   
\end{equation}
Since ${\rm low}(U_{\circ}^T M_\beta U_{\circ}) = 0$ this implies 
\begin{equation}
\label{equation konstantinov init}
{\rm low}\left( e^{-\alpha Y} U_{init}^T (\hat M_\beta - \sigma W_\beta) U_{init} e^{-\alpha Y} \right)  = 0 \qquad \Rightarrow \qquad   
{\rm low}\left( [ U^T_{init} \hat M_\beta U_{init} , \alpha Y] \right)  =  {\rm low}\left( U_{init}^T \sigma  W_\beta U_{init}\right) + O(\alpha^2)
\end{equation}
The strictly lower-diagonal part of $ [ A , \alpha Y] $ is equal to the strictly lower diagonal part of $[ A , {\rm low}(\alpha Y)]$, if $A$ is upper-triangular.
Then, by considering the projection to the subspace of strictly lower diagonal matrices of \eqref{equation konstantinov init} (see proof of Theorem \ref{theorem alpha} for more details), one obtains
\begin{equation}
T_\beta P_{\rm low}{\rm vec}(\alpha Y)  = P_{\rm low} {\rm vec}\left( U^T_{init} \sigma  W_\beta U_{init}\right)  + O(\alpha^2)
\end{equation}
with $T_\beta$ defined in \eqref{initial T}.
Since $T_\beta$ is invertible one has 
\begin{equation}
P_{\rm low}{\rm vec}(\alpha Y) = T_\beta^{-1}P_{\rm low} {\rm vec}\left( U^T_{init} \sigma  W_\beta U_{init}\right) 
\end{equation}
and taking the norm in both sides
\begin{equation}
\alpha \leq  \sqrt 2 \| T^{-1}_\beta \| \| P_{\rm low} {\rm vec}\left( U^T_{init} \sigma  \hat W_\beta U_{init}\right) \| + O(\alpha^2)
\end{equation}
where $ \| T_{\beta}^{-1}  \|_{2}$ is the spectral norm of $T^{-1}_{\beta}$. 
This implies that the initialization matrix $U_{init}$ obtained from the Schur decomposition of $\hat M_\beta$ can be written as $U_{init} = U_{\circ} e^{- \alpha Y}$, with $\alpha$ obeying
\begin{equation}
\alpha \leq \alpha_{init}  \qquad \alpha_{init}  =  \sigma \sqrt{2 N} \| T^{-1}_{\beta}  \|_{2} + O(\alpha_{init}^2)
\end{equation}
where we have used $ \| {\rm Low} {\rm vec}\left( U_{init}^T W_\beta U_{init}\right) \| \leq \sqrt N \| \beta \|  = \sqrt N$, since $\| W_n \| \leq 1$ and $\| \beta \| = 1$ by assumption.
Now, the initialization matrix $U_{init}$ belongs to the convex region containing $U_{\circ}$ if $\alpha < \alpha_{max}$, with $\alpha_{max}$  given in \eqref{alpha max}.
It follows that a descent algorithm initialized with $U_{init}$ converges to the minimum of the convex region containing $U_{\circ}$ if 
\begin{equation}
\sigma \sqrt{2 N} \| T^{-1}_{\beta}  \|_{2}  \leq \frac{2 \varepsilon - \sigma A_{\sigma}}{A_\alpha} + O(\sigma^2)
\end{equation}
or equivalently 
\begin{equation}
\sigma \leq \frac{2 \varepsilon}{ \sqrt{2N} \| T^{-1}_\beta\|_2 A_{\alpha} + A_{\sigma} } + O(\sigma^2)
\end{equation}
$\square$

\subsection{Proof of Lemma \ref{lemma observable matrices}}
Let $m_n$ be defined by $[m_n ]_{n'n''} =  {\mathbb T}_{nn'n''}$ for all $n,n',n''=1, \dots, N$.
From the definition of tensor slice $[\hat m_n ]_{n'n''} = \hat {\mathbb T}_{nn'n''}$ on has $\hat m_n = m_n + \sigma e_n$, where the noise term is defined by $[e_n]_{n'n''}=E_{nn'n''}$.
Let $m = \sum_n m_n$ and $e=\sum_n e_n$, then, from the definition of ${\mathbb T}$ given in \eqref{tensor} on has
\begin{equation}
\hat m_n = m_n + \sigma e_n\qquad m_n = Z {\rm diag}({\bf e}_n^T Z) Z^T \qquad \hat m = m + \sigma e  \qquad  m  =  Z {\rm diag}(1^T Z) Z^T  
\end{equation}
and 
\begin{equation}
\hat M_n = \hat m_n \hat m^{-1} = (m_n + \sigma e_n)(m + \sigma e)^{-1}  = m_n m^{-1} + \sigma \left( e_n m^{-1} + m_n m^{-1} e m^{-1} \right) + O(\sigma^2) 
\end{equation}
where it is easy to check that 
\begin{equation}
m_n m^{-1} =  Z {\rm diag}({\bf e}_n^T Z) Z^T \left(  Z {\rm diag}(1^T Z) Z^T\right)^{-1} = Z {\rm diag}({\bf e}_n^T Z) \left( {\rm diag}(1^T Z)\right)^{-1} Z^{-1}
\end{equation}
where we have assumed $d = N$ and the matrices $Z$ to be invertible. 
From the definitions above it follows
\begin{eqnarray}
\| m_n \|  & = & \| Z {\rm diag}({\bf e}_n^T Z) Z^T \| \\
& \leq & \| Z \|^2 \| {\rm diag}({\bf e}_n^T Z) \| \\
& \leq & \| Z \|^2 \sqrt{N} \max_n | Z_{ni} | \\
& \leq & \| Z \|^2 \sqrt{N} \max | Z |
\end{eqnarray}
and, assuming $ [1^T Z]_{i} \neq 0$ for all $i=1,\dots,d$,  
\begin{eqnarray}
\| m^{-1} \|  & = & \| \left(\sum_{n=1}^N m_n \right)^{-1} \| \\
& = & \| Z^{-T} \left( {\rm diag}(1^T Z) \right)^{-1} Z^{-1} \| \\
& \leq & \| Z^{-1} \|^2 \| ({\rm diag}(1^T Z))^{-1} \| \\
&  = & \| Z^{-1} \|^2 \frac{ \sqrt{N}}{\min | 1^T Z |} 
\end{eqnarray}  
This implies, for all $n=1,\dots,N$,  
\begin{eqnarray}
\| M_n \| & = & \| m_n m^{-1} \| \\
& \leq  & \| m_n \| \| m^{-1} \| \\
& \leq   & N \kappa(Z)^2 \frac{ \max | Z |}{\min | 1^T Z |}
\end{eqnarray} 
and
\begin{eqnarray}
\| W_n \| & = & \| e_n m^{-1} +  m_n m^{-1} e m^{-1} \|  \\
& \leq &\varepsilon \|  m^{-1} \|(1 +  \| m_n \| \| m^{-1} \|)  \\
& \leq &\varepsilon\| Z^{-1} \|^2 \frac{ \sqrt{N}}{\min | 1^T Z |} \left(1 +  N \kappa(Z)^2 \frac{ \max | Z |}{\min | 1^T Z |} \right)\\
& \leq & \varepsilon \frac{\kappa(Z)^2}{\| Z \|^2} \frac{ \sqrt{N}}{\min | 1^T Z |} \left(1 +  N \kappa(Z)^2 \frac{ \max | Z |}{\min | 1^T Z |} \right)  
\end{eqnarray}
$\square$

\subsection{Proof of Theorem \ref{main theorem}}
\label{proof main theorem}
Lemma \ref{lemma observable matrices} shows that the matrices $\hat M_n$ are approximately jointly diagonalizable.
Let ${\cal M}_{\sigma} = \{ \hat M_n \}_{n=1}^N$ and  ${\cal M}_{\circ} = \{ \hat M_n |_{\sigma=0} \}_{n=1}^N$.
Assume that ${\cal M}_{\circ}$ is such that \eqref{non-degeneracy condition lemma} is satisfied.
In this case the solutions of \eqref{joint triangularization} are characterized by the Theorem \ref{theorem alpha}.
Now, let $U_*$ be a minimizer of \eqref{joint triangularization}, then $U_*$ can be written as $U_* = U_{\circ} e^{\alpha_* X_*} $, with $\|X_*\|=1$, $X_* = -X_*^T$ and $\alpha_*$ obeying the bound given by Theorem \ref{theorem alpha}.
According to \eqref{estimation Z}, the approximate joint triangularizer $U_*$ can be used to estimate the element of the tensor component $Z$.
The distance between the estimated joint eigenvalues and the exact eigenvalues of a set of nearly jointly diagonalizable matrices is bounded by Lemma \ref{lemma eigenvalues}.     
Using the result of Theorem \ref{lemma bound non spectral} and  Lemma \ref{lemma eigenvalues} with the definition \eqref{estimation Z} one obtains
\begin{equation}
\left| \frac{Z_{ni}^*}{[1^TZ^*]_i} - \frac{Z_{ni}}{[1^TZ]_i} \right| \leq 4 N \sigma  \frac{\sqrt{d(d-1)} \kappa(V)^4}{\gamma} {\rm M}^2 {\rm W} + \sigma {\rm W} + O(\sigma^2)
\end{equation}
for all $i=1, \dots, d$ and all $n=1, \dots,N$.
From Lemma \ref{lemma observable matrices} on has 
\begin{equation}
{\rm M}  \leq N \kappa(Z)^2 \frac{ \max | Z |}{\min | 1^T Z |} \qquad
{\rm W}  \leq   \varepsilon \frac{\kappa(Z)^2}{\| Z \|^2} \frac{ \sqrt{N}}{\min | 1^T Z |} \left(1 +  N \kappa(Z)^2 \frac{ \max | Z |}{\min | 1^T Z |} \right) 
\end{equation}
from which the claim of the theorem.
$\square$ 

\subsection{Auxiliary lemmas}
\begin{lemma}
\label{lemma linear combination}
If \eqref{non-degeneracy condition lemma} holds it is possible to find $\beta = [\beta_1, \dots, \beta_N]$ such that 
\begin{equation}
M = \sum_{n=1}^N \beta_n M_n
\end{equation}
has real distinct eigenvalues.
\end{lemma}
\paragraph{Proof of Lemma \ref{lemma linear combination}}
Let $\beta = [\beta_1, \dots, \beta_N]$, then the eigenvalues of $M = \sum_{n=1}^N \beta_n M_n$ are
\begin{equation}
\lambda_i(M) = \sum_{n=1}^N \beta_n \Lambda_{in}  \qquad i=1, \dots, d
\end{equation} 
We want to show that \eqref{non-degeneracy condition lemma} implies that it is possible to find $\beta_1, \dots, \beta_N$ such that $\lambda_i \neq \lambda_{i'}$ for all $i\neq i'$, with $i,i'=1,\dots,d$.
This can be seen as follows. 
It is aways possible to choose $\tilde m_2$ such that $\lambda_1(\tilde m_2) \neq \lambda_2(\tilde m_2)$.
Now, assume that $\tilde m_n$ is such that $\lambda_i(\tilde m_n) \neq \lambda_j(\tilde m_n)$ for all $i\neq j$ and $i,j \leq n$.
Consider $\lambda_{n+1}(\tilde m_n)$.
We want to show that it is possible to find a matrix $m_{n+1}$ and a coefficient $\beta_{n+1}$ such that the first $n+1$ eigenvalues of  $\tilde m_{n+1} = \tilde m_n + \beta_{n+1} m_{n+1}$ are distinct, that is $\lambda_i(\tilde m_{n+1}) \neq \lambda_j(\tilde m_{n+1})$ for all $i\neq j$ and $i,j \leq n+1$.
If $\lambda_{n+1}(\tilde m_n) \neq \lambda_i(\tilde m_n)$ for all $i \leq n$, one has $\tilde m_{n+1} = \tilde m_n$.
Otherwise, there exists an $i \leq n$ such that $\lambda_{n+1}(\tilde m_n)  =  \lambda_i(\tilde m_n)$.
Note that, since $\lambda_i(\tilde m_n) \neq \lambda_j(\tilde m_n)$ for all $i\neq j$ and $i,j \leq n$, there is only one such $i$.
Let $m_{n+1}$ be the matrix in ${\cal M}_{\circ}$ satisfying $\lambda_{n+1}(\tilde m_n)  \neq \lambda_i(\tilde m_n)$ and 
\begin{equation}
\beta_{n+1} \in {\mathbf R} \qquad {\rm s.t.} \qquad \beta_{n+1} \neq 0 \quad {\rm and} \quad \beta_{n+1} \neq \frac{\lambda_{i}(\tilde m_{n}) -  \lambda_{j}(\tilde m_{n})}{\lambda_{j}(m_{n+1}) -  \lambda_{i}(m_{n+1})} \qquad {\rm for \ all} \quad i \neq j \quad i,j \leq n
\end{equation}
Then it is easy to check that the first $n+1$ eigenvalues of $\tilde m_{n+1} = \tilde m_{n} + \beta_{n+1} m_{n+1}$ are distinct. 
The matrix $M$ is then constructed by repeating the above procedure until $n+1 = d$.
$\square$

%
\begin{lemma}
\label{lemma stationarity equation}
Let $U$ be a stationary point of \eqref{joint triangularization}, then
\begin{equation}
\label{stationarity equation lemma}
S - S^T   = 0 \qquad S = \sum_{n=1}^N \left[ U^T \hat M^T_n U , {\rm low}( U^T \hat M_n U)  \right]
\end{equation}
\end{lemma}  

\paragraph{Proof of Lemma \ref{lemma stationarity equation}}
Let $f(U)$ be a function defined on ${\mathbb O}(d)$.
The directional derivatives of $f$ at $U$ in the direction $X$ are defined as
\begin{eqnarray}
 D_X f(U) &=& \langle \nabla f(U), X \rangle \\
& = & \frac{d}{dt}f(U e^{Xt})|_{t=0}
 \end{eqnarray}  
where $X = -X^T$ and the scalar product in the tangent space is defined by $\langle A,B \rangle= {\rm Tr}(A^T B)$.
In particular, for \eqref{joint triangularization} one has
\begin{eqnarray}
\langle X , \nabla {\cal L}(U,{\cal M}_{\sigma}) \rangle  & = &\left. \frac{d}{dt}{\cal L}(U e^{Xt},{\cal M}_{\sigma}) \right|_{t=0}  \\
& = & \sum_{n=1}^N {\rm Tr} \left( [U^T \hat M^T_n U,X]{\rm low}(U^T \hat M_n U) +  [U^T \hat M_n U,X]{\rm up}(U^T \hat M^T_n U)  \right) \\ 
& = -  & \sum_{n=1}^N {\rm Tr} \left( X[ U^T \hat M^T_n U,{\rm low}(U^T \hat M_n U)] + 
X[U^T \hat M_n U,{\rm up}(U^T \hat M^T_n U)]  \right) \\ 
\label{gradient identification}
& = & \langle X , S - S^T \rangle
\end{eqnarray} 
where we have defined $S = \sum_{n=1}^N [ U^T \hat M^T_n U,{\rm low}(U^T \hat M_n U)]$ and used $\sum_{n=1}^N [U^T \hat M_n U,{\rm up}(U^T \hat M^T_n U)] = - S^T $.
From \eqref{gradient identification} one has $\nabla {\cal L}(U,{\cal M}_{\sigma}) = S - S^T$ and \eqref{stationarity equation lemma} follows from the stationarity condition $\nabla {\cal L} = 0$. 
$\square$\\
\begin{lemma}
\label{lemma T}
Let ${\cal M}_{\circ} = \{ M_n  = V {\rm diag}([\Lambda_{n1}, \dots , \Lambda_{nd}]) V^{-1} \}_{n=1}^N$ be a set of jointly diagonalizable matrices such that 
\begin{equation}
\gamma = \min_{i>i'} \sum_{n=1}^N (\Lambda_{ni} - \Lambda_{ni'})^2 >0
\end{equation}
and let $U_{\circ}$ be an exact triangularizer of ${\cal M}_{\circ} $.
Then the operator 
\begin{equation}
T = \sum_n t^T_n t_n  \qquad     t_n = P_{\rm low} (1 \otimes U_{\circ}^T M^T_n U_{\circ} - U_{\circ}^T M_n U_{\circ} \otimes 1 )P_{\rm low}^T
\end{equation}
is invertible and 
\begin{equation}
\sigma_{\min}(T) \geq  \frac{\gamma}{\kappa(V)^4} \qquad  \| T^{-1} \| \leq   \sqrt{\frac{d(d-1)}{2}}  \frac{\kappa(V)^4}{\gamma} 
\end{equation} 
\end{lemma}
\paragraph{Proof of Lemma \ref{lemma T}}
Since $U_{\circ}^T M_n U_{\circ}$ is upper-triangular for all $n=1, \dots, N$, the matrices $(1 \otimes U_{\circ}^T M^T_n U_{\circ} - U_{\circ}^T M_n U_{\circ} \otimes 1)$ are block upper-triangular matrices and their diagonal blocks are lower triangular.
For all $n=1, \dots, N$ one has $U_{\circ}^T M_n U_{\circ} = U_{\circ}^T V \Lambda_n V^{-1} U_{\circ}$ where we have defined $\Lambda_n =  {\rm diag}([\Lambda_{n1}, \dots, \Lambda_{nd}])$. 
Then 
\begin{eqnarray}
t_n &=& P_{\rm low} (1\otimes U_{\circ}^T M^T_n U_{\circ} - U_{\circ}^T M_n U_{\circ} \otimes 1) P_{\rm low} \\
& = & P_{\rm low} (U_{\circ}^T V \otimes U_{\circ}^TV^{-T}) (1 \otimes  \Lambda_n -  \Lambda_n \otimes 1) (V^{-1}U_{\circ}\otimes V^{T}U_{\circ}) P_{\rm low}^T  \\
& = & P_{\rm low} (U_{\circ}^T V \otimes U_{\circ}^TV^{-T})  P_{\rm low}^T P_{\rm low} (1 \otimes  \Lambda_n -  \Lambda_n \otimes 1) P_{\rm low}^T P_{\rm low} (V^{-1}U_{\circ}\otimes V^{T}U_{\circ}) P_{\rm low}^T   \\
& = & \tilde V \Gamma_n \tilde V^{-1}
\end{eqnarray}
where we have defined $\Gamma_n =P_{\rm low} (1 \otimes  \Lambda_n -  \Lambda_n \otimes 1) P_{\rm low}^T $ and $\tilde V  = P_{\rm low} (U_{\circ}^TV \otimes U_{\circ}^TV^{-T})  P_{\rm low}^T$, $\tilde V^{-1}  = P_{\rm low}  (V^{-1}U_{\circ}\otimes V^{T}U_{\circ}) P_{\rm low}^T $ and the last equality follows form the fact that $U_{\circ}^T V$ is upper triangular (see Lemma \ref{lemma low insertion}).
The positive semi-definite matrix $T$ can be rewritten as
\begin{equation}
T = \sum_{n=1}^N t_n^T t_n = W^T W \qquad W = [\tilde V^{-T} \Gamma_1 \tilde V^T , \dots, \tilde V^{-T} \Gamma_N \tilde V^T]^T  = (1 \otimes \tilde V) [\Gamma_1 , \dots, \Gamma_N ]^T \tilde V^{-1} 
\end{equation}
A bound on the smallest singular value of $T$ can be obtained as follows 
\begin{eqnarray}
\label{sigma min T}
\sigma_{min}(T) &=&  \sigma_{min}\left( \tilde V^{-T} [\Gamma_1, \dots , \Gamma_n] (1 \otimes \tilde V^{T}) (1 \otimes \tilde V) [\Gamma_1, \dots , \Gamma_n]^T \tilde V^{-1} \right) \\
 &\geq & \sigma_{min}(\tilde V^{-1})^2  \sigma_{min}\left( [\Gamma_1, \dots , \Gamma_n] (1 \otimes \tilde V^{T}) (1 \otimes \tilde V) [\Gamma_1, \dots , \Gamma_n]^T) \right) \\
 & =  & \sigma_{min}(\tilde V^{-1})^2  \left( \min_{\|x \| = 1} x^T [\Gamma_1, \dots , \Gamma_n] (1 \otimes \tilde V^{T}) (1 \otimes \tilde V)  [\Gamma_1, \dots , \Gamma_n]^T x \right) \\
 & \geq   & \sigma_{min}(\tilde V^{-1})^2 \sigma_{min}(1 \otimes \tilde V)^2  \left( \min_{\|x \| = 1} x^T [\Gamma_1, \dots , \Gamma_n]  [\Gamma_1, \dots , \Gamma_n]^T x \right) \\
 & \geq   & \sigma_{min}(\tilde V^{-1})^2 \sigma_{min}(\tilde V)^2  \left( \min_{\|x \| = 1} x^T [\Gamma_1, \dots , \Gamma_n]  [\Gamma_1, \dots , \Gamma_n]^T x \right) \\
 & \geq   & \sigma_{min}(\tilde V^{-1})^2 \sigma_{min}(\tilde V)^2  \left( \min_{\|x \| = 1} x^T  {\rm diag}\left( \left[ \sum_n [\Gamma_n]^2_{11}, \dots, \sum_n [\Gamma_n]^2_{\tilde d \tilde d } \right] \right)  x \right) 
 \end{eqnarray}
 where we have defined $\tilde d = \frac{d(d-1)}{2} $.
The minimization problem between brackets is solved by ${\bf e}_{i_*} $ with $i_* = {\rm arg} \min_{i} \sum_{n=1}^N [\Gamma_n]^2_{ii} $ and one has 
\begin{equation}
\gamma = {\bf e}_{i_*}^T  {\rm diag}\left( \left[ \sum_n [\Gamma_n]^2_{11}, \dots, \sum_n [\Gamma_n]^2_{\tilde d \tilde d } \right] \right)  {\bf e}_{i_*} =  \min_{j<j'} \sum_{n=1}^N (\Lambda_{nj} - \Lambda_{nj'})^2 
\end{equation}
where $i_*$ and $(j_*,j_*') = {\rm arg} \min_{j<j'} \sum_{n=1}^N (\Lambda_{nj} - \Lambda_{nj'})^2$ are related by
\begin{equation}
 i = f(j,j') \qquad f(j,j') = \sum_{k=1}^{j-1}(d-k)+j'-j  \qquad {\rm for} \quad j < j'
\end{equation}
This implies
\begin{equation}
\sigma_{\min}(T) \geq \frac{\gamma}{\kappa(V)^4} \qquad 
\|T^{-1} \| \leq   \sqrt{\frac{d(d-1)}{2} } \frac{\kappa(V)^4}{ \gamma}   
\end{equation} 
where we have used 
\begin{equation}
\sigma_{min}(\tilde V)  = \sigma_{min}(P_{\rm low}(V\otimes V^{-T}) P_{\rm low}^T ) \geq \sigma_{min}(V) \sigma_{min}(V^{-1})  = \frac{ \sigma_{min}(V)}{\sigma_{max}(V)} = \frac{1}{\kappa(V)} 
\end{equation}
\begin{equation}
\sigma_{min}(\tilde V^{-1}) = \sigma_{min}(P_{\rm low}(V^{-1}\otimes V^T) P_{\rm low}^T )  \geq \sigma_{min}(V^{-1}) \sigma_{min}(V)  =  \frac{\sigma_{min}(V)}{\sigma_{max}(V)} = \frac{1}{\kappa(V)} 
\end{equation}
$\square$ \\

\begin{lemma}
\label{lemma low insertion}
Let $A$ be an upper triangular (invertible) matrix and $\Sigma$ a diagonal matrix, then for any $B$   
\begin{equation}
{\rm Low} (A \otimes A^{-T}) (1 \otimes \Sigma -\Sigma \otimes 1)   (A^{-1} \otimes A^{T}) {\rm Low}  = {\rm Low} (A \otimes A^{-T}) {\rm Low}  (1 \otimes \Sigma -\Sigma \otimes 1){\rm Low}  (A^{-1} \otimes A^{T}) {\rm Low}   
 \end{equation}
\end{lemma}

\paragraph{Proof of Lemma \ref{lemma low insertion}}
Let $B$ be any matrix of the same dimension as $A$,
\begin{eqnarray}
M_1 & = & {\rm mat}\left( {\rm Low} (A \otimes A^{-T})  (1 \otimes \Sigma -\Sigma \otimes 1) (A^{-1} \otimes A^{T}) {\rm Low} \right) \\
& = &{\rm mat}\left( {\rm Low}  (1 \otimes \Sigma -\Sigma \otimes 1) (A^{-1} \otimes A^{T}){\rm vec}\left( {\rm low}(B) \right) \right) \\
& = &{\rm mat}\left( {\rm Low} (A \otimes A^{-T}) (1 \otimes \Sigma -\Sigma \otimes 1) {\rm vec}\left( A^T {\rm low}(B) A^{-T} \right) \right) \\
& = &{\rm mat}\left( {\rm Low} (A \otimes A^{-T}) {\rm vec}\left( [\Sigma,  A^T {\rm low}(B) A^{-T}] \right) \right)\\
& = &{\rm mat}\left( {\rm Low} {\rm vec}\left( A^{-T}  [\Sigma,  A^T {\rm low}(B) A^{-T} ] A^{T}  \right) \right) \\
& = &{\rm low}\left( A^{-T}  [\Sigma, A^{T} {\rm low}(B) A^{-T}]  A^{T}  \right) 
\end{eqnarray} 
and 
\begin{eqnarray}
M_2 & = & {\rm mat}\left( {\rm Low} (A \otimes A^{-T}) {\rm Low}(1 \otimes \Sigma -\Sigma \otimes 1){\rm Low} (A^{-1} \otimes A^{T}){\rm Low} {\rm vec}(B) \right) \\
& = &{\rm mat}\left( {\rm Low} (A \otimes A^{-T}) {\rm Low}  (1 \otimes \Sigma -\Sigma \otimes 1) {\rm vec}\left( {\rm low}\left( A^{T} {\rm low}(B) A^{-T}\right) \right) \right) \\
& = &{\rm mat}\left( {\rm Low} (A \otimes A^{-T}) {\rm Low} {\rm vec}\left( [\Sigma, {\rm low}\left( A^{T} {\rm low}(B) A^{-T}\right) ] \right) \right)\\
& = &{\rm mat}\left( {\rm Low} {\rm vec}\left(A^{-T}{\rm low}\left( [\Sigma, {\rm low}\left(A^{T} {\rm low}(B) A^{-T} \right) ]\right) A^{T} \right)  \right)\\
& = & {\rm low}\left( A^{-T}  \ {\rm low}\left([\Sigma, {\rm low}\left(A^{T} {\rm low}(B) A^{-T}\right) ] \right)  A^{T}  \right) 
\end{eqnarray} 
Then $M_1 = M_2$ can be shown by observing that  $A^{T} {\rm low}(B) A^{-T}$ is a lower-diagonal matrix if $A$ is upper triangular. 
Then, for every lower-diagonal matrix $C$, one has 
\begin{equation}
[\Sigma, C] = [\Sigma, {\rm low}(C) + {\rm diag}(C)] = [\Sigma, {\rm low}(C) ]  = {\rm low}([\Sigma, {\rm low}(C) ])
\end{equation}
because diagonal matrices always commute and the commutator of a strictly lower diagonal matrix with a diagonal matrix is strictly lower diagonal. $\square$\\

\begin{lemma}
\label{lemma T beta}
Let $A$ be an upper triangular matrix with real nonzero eigenvalues.
If  $A$ is invertible and the eigenvalues of $A$ satisfy $\lambda_{i}(A) \neq \lambda_{i'}(A)$ for all $i \neq i'$ the matrix
\begin{equation}
T_{A} = P_{\rm low} (1 \otimes A^T - A \otimes 1) P_{\rm low}^T
\end{equation}
is invertible.
\end{lemma}

\paragraph{Proof of Lemma \ref{lemma T beta}}
From the spectral decomposition of the matrix $A$ one has $A = V \Lambda V^{-1} $, with $V$ upper triangular and $\Lambda$ diagonal, and 
\begin{eqnarray}
T_{A} &=& P_{\rm low} (V \otimes V^{-T})(1 \otimes \Lambda - \Lambda \otimes 1)(V^{-1} \otimes V^{T}) P_{\rm low}^T \\
&=& P_{\rm low} (V \otimes V^{-T}){\rm Low}(1 \otimes \Lambda -\Lambda \otimes 1){\rm Low} (V^{-1} \otimes V^{T}) P_{\rm low}^T 
\end{eqnarray}
where the second equality follows from the fact that $ (V^{-1} \otimes V^{T}) P_{\rm low}^T \tilde a = {\rm Low}(V^{-1} \otimes V^{T}) P_{\rm low}^T \tilde a$ for any $\frac{d(d-1)}{2}$-dimensional vector $\tilde a$ and $(1 \otimes \Lambda -\Lambda \otimes 1){\rm Low} a  =  {\rm Low} (1 \otimes \Lambda -\Lambda \otimes 1){\rm Low} a $ for any $d$-dimensional vector $a$.
The smallest singular value of $T_A$ obeys  
\begin{eqnarray}
\sigma_{\min}(T_A)  &\geq & C_1^2 C_2
\end{eqnarray}
where
\begin{eqnarray}
C_1  = \sigma_{\min}(V^{-1})\sigma_{\min}(V) = \frac{\sigma_{\min}(V)}{\sigma_{\max}(V)}
\qquad
C_2  = \min \{ \| P_{\rm low}  (1 \otimes \Lambda -\Lambda \otimes 1)P_{\rm low}^T x \| , \|x \|=1\}  =   \min_{i<i'} | \lambda(A)_{i}-\lambda(A)_{i'} |
\end{eqnarray}
This implies that $T_A$ is invertible if $V$ is full rank and $\lambda_i(A) \neq \lambda_{i'}(A)$ for all $i \neq i'$, which are both fulfilled by assumption.
$\square$ 

\begin{lemma}
\label{lemma hessian}
The Hessian of ${\cal L}$ at $U = U_{\circ} e^{\alpha Y}$, where  $U_{\circ}$ is an exact triangularizer of ${\cal M}_{\circ} = \{\hat M_n |_{\sigma=0} \}_{n=1}^N$ and $Y = -Y^T$, $\| Y \| =1$,   is positive definite for all $Y$ if 
\begin{equation}
\alpha \leq \alpha_{max}   \qquad \alpha_{max} = \frac{2 \varepsilon - \sigma A_{\sigma}}{A_\alpha} + O((\alpha+\sigma)^2)
\end{equation} 
\begin{equation}
\varepsilon = \frac{\gamma}{2 \kappa(V)^4}  \qquad \gamma = \min_{j<j'} \sum_{n=1}^N(\Lambda_{nj}-\Lambda_{nj'})^2
\qquad 
A_\alpha = 32 \sum_{n=1}^N \| M_n\|^2  \qquad 
A_\sigma = 16  \sqrt{N} \sqrt{\sum_{n=1}^N \| M_n \|^2}
\end{equation}
\end{lemma}

\paragraph{Proof of Lemma \ref{lemma hessian}}
Let 
\begin{equation}
{\cal L}(U,{\cal M_{\sigma}}) = \sum_{n=1}^N {\rm Tr}(g_n^T g_n) \qquad g_n= {\rm low}(U^T \hat M_n U)
\end{equation}
Then we have $\langle X , \nabla {\cal L}(U) \rangle = \frac{d}{dt} {\cal L}(Ue^{t X})|_{t=0} = \sum_{n=1}^N {\rm Tr}(\dot g_n^T g_n + g_n^T \dot g_n)$ where $X=-X^T$ and $\dot g_n= \frac{d}{dt} g_n(Ue^{tX}) |_{t=0} = {\rm low}([U^T \hat M_n U,X])$.
The second derivative in the direction $X$ defines the Hessian of ${\cal L}$ at $U$ via 
\begin{equation}
\label{second order directional derivative}
\langle X , \nabla^2 {\cal L}X \rangle  = \frac{d^2}{dt^2} {\cal L}(Ue^{t X})|_{t=0}  = \sum_{n=1}^N {\rm Tr}( 2 \dot g_n^T \dot g_n + \ddot g_n^T g_n + g_n^T \ddot g_n)
\end{equation}
where $\ddot g_n = \frac{d^2}{dt^2} g_n(Ue^{tX})|_{t=0} = {\rm low}([[U^T \hat M_n U,X],X])$.
Let $f(U, {\cal M}_{\sigma})$ be a general function of $U = U_{\circ} e^{\alpha Y}$, where  $U_{\circ}$ is an exact triangularizer of ${\cal M}_{\circ} = \{\hat M_n |_{\sigma=0} \}_{n=1}^N$ and $Y = -Y^T$, $\| Y \| =1$, and the empirical matrices $\hat M_n$.
The double expansion, respect to the parameter $\alpha$ and $\sigma$ is
\begin{equation}
f = f |_{(\alpha=0, \sigma=0)} + \alpha \partial_\alpha f |_{\sigma=0} + \sigma \partial_\sigma f |_{\alpha=0} + O((\alpha + \sigma)^2)
\end{equation}
Now, consider the double expansion of the functions $g_n$, $\dot g_n$ and $\ddot g_n$.
In the first order approximation one obtains
\begin{eqnarray}
\label{expansion hessian}
\langle X , \nabla^2 {\cal L} X \rangle  &=& \sum_{n=1}^N {\rm Tr} \left(
 2 \dot g_n^T\dot g_n  + \sigma (\dot g_n^T \partial_{\sigma} \dot g_n + \partial_{\sigma}\dot g_n^T \dot g_n  ) + \alpha ( \dot g_n^T \partial_{\alpha} \dot g_n + \dot g_n^T\partial_{\alpha}\dot g_n) +\right.\\
 && \left. \sigma(\ddot g_n^T \partial_\sigma g_n + \partial_\sigma g_n^T \ddot g_n ) + \nonumber
   \alpha(\ddot g_n^T \partial_{\alpha} g_n + \partial_\alpha g_n^T \ddot g_n ) 
   \right) + O((\alpha + \sigma)^2)   
\end{eqnarray}
where the first term is always nonnegative.
Now, the Hessian of ${\cal L}$ at $U$ is positive definite if $\langle X , \nabla^2 {\cal L} X \rangle$, for all $X$ such that $X = -X^T$.
The non negativity of \eqref{expansion hessian} is guaranteed by the following condition 
\begin{equation}
2 \sum_{n=1}^N {\rm Tr} (\dot g_n^T\dot g_n ) \geq \alpha  \tilde A_\alpha + \sigma \tilde A_{\sigma} +O((\alpha+\sigma)^2)
\end{equation}
where
\begin{eqnarray}
  \tilde A_\alpha  =  \left| \sum_{n=1}^N {\rm Tr} \left( \dot g_n^T \partial_{\alpha} \dot g_n+ \dot g_n^T\partial_{\alpha}\dot g_n + \ddot g_n^T \partial_{\alpha} g_n+  \partial_\alpha g_n^T \ddot g_n \right) \right|\qquad   
\tilde A_{\sigma} =  \left| \sum_{n=1}^N {\rm Tr}  \left( \dot g_n^T \partial_{\sigma} \dot g_n +\partial_{\sigma}\dot g_n^T \dot g_n + \ddot g_n^T \partial_\sigma g_n+ \partial_\sigma g_n^T \ddot g_n \right) \right| 
\end{eqnarray}
We seek some $\varepsilon$, $A_\alpha$ and  $A_\sigma$ such that 
\begin{equation}
\label{condition upper-lower bounds}
\sum_{n=1}^N {\rm Tr} (\dot g_n^T\dot g_n ) \geq  \varepsilon \| X \|^2  \qquad A_\alpha\| X \|^2    \geq \tilde A_\alpha \qquad  A_{\sigma} \| X \|^2  \geq  \tilde A_{\sigma}
\end{equation}
Given $\varepsilon$, $A_\alpha$ and  $A_\sigma$ satisfying \eqref{condition upper-lower bounds}, the non negativity of the Hessian is implied by 
\begin{equation}
2 \varepsilon \geq \alpha A_{\alpha} + \sigma A_{\sigma}
\end{equation}
from which the condition on $\alpha$ stated by the lemma. 
The explicit form of $\varepsilon$, $A_\alpha$ and  $A_\sigma$ are provided by Lemma \ref{lemma varepsilon} and  Lemma \ref{lemma A upper bounds}.
$\square$
\begin{lemma}
\label{lemma varepsilon}
A possible choice of $\varepsilon >0$ satisfying \eqref{condition upper-lower bounds} is given by
\begin{equation}
 \varepsilon = \frac{\gamma}{2 \kappa(V)^4}  \qquad \gamma =  \min_{j<j'} \sum_{n=1}^N(\Lambda_{nj}-\Lambda_{nj'})^2   
\end{equation}
with $V$ and $\Lambda$ defined in \eqref{jointschur}.
\end{lemma}
\paragraph{Proof of Lemma \ref{lemma varepsilon}}
This can be seen as follows:
\begin{eqnarray}
\sum_{n=1}^N {\rm Tr}( \dot g_n^T\dot g_n) & = & \sum_{n=1}^N {\rm Tr}\left(  {\rm low}\left([U_{\circ}^T M_n U_{\circ},X]\right)^T  {\rm low}\left([U_{\circ}^T M_n U_{\circ},X] \right) \right)\\
& = & \sum_{n=1}^N {\rm Tr}\left( {\rm low}\left([U_{\circ}^T M_n U_{\circ},{\rm low}(X)]\right)^T  {\rm low}\left([U_{\circ}^T M_n U_{\circ},{\rm low}(X)] \right) \right)\\
& = & \sum_{n=1}^N {\rm vec}\left( {\rm low}\left([U_{\circ}^T M_n U_{\circ},{\rm low}(X)]\right) \right)^T {\rm vec}\left( {\rm low}\left([U_{\circ}^T M_n U_{\circ},{\rm low}(X)]\right) \right) \\
& = & \sum_{n=1}^N {\rm vec}(X)^T {\rm Low} (1 \otimes U_{\circ}^T M^T_n U_{\circ} - U_{\circ}^T M^T_n U_{\circ}\otimes 1)^T {\rm Low}\\
&&\quad (1 \otimes U_{\circ}^T M^T_n U_{\circ} - U_{\circ}^T M^T_n U_{\circ}\otimes 1) {\rm Low} {\rm vec}(X) \qquad \\
& = & {\rm vec}(X)^T P_{\rm low} \left(\sum_{n=1}^N t^T_n t_n\right) P_{\rm low}^T {\rm vec}(X)\\
& = & {\rm vec}(X)^T P_{\rm low} T  P_{\rm low}^T {\rm vec}(X)
\end{eqnarray}
where we have used ${\rm Low} = P_{\rm low} P_{\rm low}^T$ and the definition of $T$ given in Lemma \ref{lemma T}.
For every $X$ such that  $X = -X^T$ one has $\| {\rm low}(X)\|  = \frac{1}{\sqrt 2} \| X \|$.
In particular
\begin{equation}
{\rm vec}(X)^T P_{\rm low} T  P_{\rm low} {\rm vec}(X) \geq \frac{1}{2} \| X \|^2 \sigma_{\min}(T) 
\end{equation} 
and using the result of Lemma \ref{lemma T} one obtains
\begin{equation}
\sum_{n=1}^N {\rm Tr}( \dot g_n^T\dot g_n)  \geq \frac{\gamma}{2 \kappa(V)^4} \| X \|^2
\end{equation}
and hence $\varepsilon = \frac{\gamma}{2 \kappa(V)^4} $. $\square$
\begin{lemma}
\label{lemma A upper bounds}
A possible choice of $A_\alpha$ and $A_{\sigma}$ satisfying \eqref{condition upper-lower bounds} is given by 
\begin{equation}
A_\alpha = 32  \sum_{n=1}^N \| M_n\|^2  \qquad 
A_\sigma = 16  \sqrt{N} \sqrt{\sum_{n=1}^N \| M_n \|^2}
\end{equation}
\end{lemma}
\paragraph{Proof of Lemma \ref{lemma A upper bounds}}
Let $a_\alpha$, $b_\alpha$, $a_\sigma$ and $b_\sigma$ be defined by 
\begin{eqnarray}
 \sum_{n=1}^N {\rm Tr} (\dot g_n^T \partial_{\alpha} \dot g_n) 
& = & \sum_{n=1}^N {\rm Tr}\left( {\rm low}([U_{\circ}^T M_n U_{\circ} ,X])^T {\rm low}([[U_{\circ}^T M_n U_{\circ} ,Y],X]) \right) \\
& \leq & \sqrt{ \sum_{n=1}^N \| {\rm low}([U_{\circ}^T M_n U_{\circ} ,X]) \|^2} \sqrt{ \sum_{n=1}^N \| {\rm low}([[U_{\circ}^T M_n U_{\circ} ,Y],X]\|^2 }\\
& \leq & \| X \|^2 \sqrt{ \sum_{n=1}^N 4 \| M_n \|^2} \sqrt{ \sum_{n=1}^N 16 \| M_n \|^2 } \\
& \leq & 8 \| X \|^2  \sum_{n=1}^N \| M_n \|^2\\
& = & a_\alpha \| X \|^2
\end{eqnarray}
\begin{eqnarray}
\sum_{n=1}^N {\rm Tr} (\ddot g_n^T \partial_{\alpha} g_n) 
 & = &\sum_{n=1}^N {\rm Tr}\left( {\rm low}([[U_{\circ}^T M_n U_{\circ},X],X])^T {\rm low}([U_{\circ}^T M_n U_{\circ} ,Y]) \right) \\
& \leq & 8 \| X \|^2  \sum_{n=1}^N \| M_n \|^2\\
& = &b_\alpha \| X \|^2
\end{eqnarray}
\begin{eqnarray}
\sum_{n=1}^N {\rm Tr} (\dot g_n^T \partial_{\sigma} \dot g_n) 
& = & \sum_{n=1}^N {\rm Tr}\left( {\rm low}([U_{\circ}^T M_n U_{\circ} ,X])^T {\rm low}([U_{\circ}^T W_n U_{\circ},X]) \right) \\
& \leq & \sqrt{ \sum_{n=1}^N \| {\rm low}([U_{\circ}^T M_n U_{\circ} ,X]) \|^2} \sqrt{ \sum_{n=1}^N \| {\rm low}([U_{\circ}^T W_n U_{\circ} ,X]\|^2 }\\
& \leq & \| X \|^2 \sqrt{ \sum_{n=1}^N 4 \| M_n \|^2} \sqrt{ \sum_{n=1}^N 4 \| W_n \|^2 } \\
& \leq & 4 \| X \|^2 \sqrt{N} \sqrt{\sum_{n=1}^N \| M_n \|^2}\\
 & = &a_{\sigma} \| X \|^2
\end{eqnarray}
\begin{eqnarray}
\sum_{n=1}^N {\rm Tr} (\ddot g_n^T \partial_{\sigma} g_n) 
 & = &\sum_{n=1}^N {\rm Tr}\left( {\rm low}([[U_{\circ}^T M_n U_{\circ} ,X],X])^T {\rm low}(U_{\circ}^T W_n U_{\circ}) \right) \\
& \leq & 4 \| X \|^2 \sqrt{N} \sqrt{\sum_{n=1}^N \| M_n \|^2}\\
 & = & b_{\sigma} \| X \|^2
\end{eqnarray}
where we have defined $a_\alpha = 8 \sum_{n=1}^N \| M_n \|^2 = b_\alpha$,   $a_\sigma = 4  \sqrt{N} \sqrt{\sum_{n=1}^N \| M_n \|^2} = b_\sigma$, used $\| Y \| = 1$ and 
\begin{eqnarray}
\sum_{n=1}^{N} {\rm Tr}(A_n B_n) &=& \sum_{n=1}^{N} {\rm vec}(A^T_n)^T {\rm vec}(B_n) \\
 &=& {\rm Tr} \left(  [{\rm vec}(A^T_1),\dots, {\rm vec}(A^T_N)]^T [ {\rm vec}(B_1), \dots, {\rm vec}(B_N)] \right) \\
 &=& {\rm vec}\left(  [{\rm vec}(A^T_1),\dots, {\rm vec}(A^T_N)]\right)^T  {\rm vec}\left([ {\rm vec}(B_1), \dots, {\rm vec}(B_N)] \right) \\
 &\leq & \| {\rm vec}\left(  [{\rm vec}(A^T_1),\dots, {\rm vec}(A^T_N)]\right) \| \|  {\rm vec}\left([ {\rm vec}(B_1), \dots, {\rm vec}(B_N)] \right) \| \\
 & = & \sqrt{\sum_{n=1}^{N} \| {\rm vec}(A^T_n) \|^2 }\sqrt{\sum_{n=1}^{N}  \| {\rm vec}(B_n) \|^2} \\
 & = & \sqrt{\sum_{n=1}^{N} \| A_n \|^2 } \sqrt{\sum_{n=1}^{N}  \|B_n \|^2} 
 \end{eqnarray}
Then we have 
\begin{equation}
\tilde A_\alpha  \leq 2 \| X \|^2(  a_\alpha  +  b_\alpha ) \qquad \tilde A_\sigma  \leq 2 \| X \|^2 ( a_\sigma  +  b_\sigma )
\end{equation}
$\square$
\begin{lemma}
\label{lemma eigenvalues}
Let $U$ and $U_{\circ}$ be respectively the approximate joint triangularizers of ${\cal M}_{\sigma}$ and the exact joint triangularizer of ${\cal M}_{\circ}$ defined in Theorem \ref{theorem alpha}.
For all $n=1, \dots, N$ and all $i=1, \dots,d $, let $\hat \lambda_i(\hat M_n) = [U^T \hat M_n U]_{ii}$ and $\lambda_i(M_n) = [U^T_{\circ} M_n U_{\circ}]_{ii}$.
Then, for all $n=1, \dots, N$ and all $i=1, \dots, d$ ,
\begin{equation}
\label{bound eigenvalues lemma}
\left| \hat \lambda_i(\hat M_n)  - \lambda_i(M_n)\right| \leq  2 \alpha \| M_n \| + \sigma \| W_n \| + O(\alpha^2)
\end{equation}
with $\alpha$ defined in Theorem \ref{theorem alpha}.
\end{lemma}

\paragraph{Proof of Lemma \ref{lemma eigenvalues}}
 Let $U$ and $U_{\circ}$ be respectively the approximate joint triangularizers of ${\cal M}_{\sigma}$ and the exact joint triangularizer of ${\cal M}_{\circ}$ defined in Theorem \ref{theorem alpha}.
Then $ U = U_{\circ} e^{\alpha X} $ with $X = -X^T$, $ \|X \|  = 1 $ and $\alpha >0$ obeying \eqref{bound alpha theorem}.
Neglecting all second order terms one has 
\begin{eqnarray}
\left| \hat \lambda_i(\hat M_n)  - \lambda_i(M_n)\right|  &=& \left| [U^T \hat M_n U]_{ii} - [U_{\circ}^T  M_n U_{\circ}]_{ii} \right|   \\
&=& \left| [e^{-\alpha X} U_{\circ}^T (M_n+\sigma W_n) U_{\circ} e^{-\alpha X} ]_{ii} - [U_{\circ}^T  M_n U_{\circ}]_{ii} \right|
\\
&=& \left| [  U^T_{\circ} M_n U_{\circ} \alpha X - \alpha X U^T_{\circ} M_n U_{\circ} ]_{ii} +\sigma [U_{\circ}^T W_n U_{\circ}]_{ii} \right| + O(\alpha^2) \\
& \leq &  2 \alpha \| M_n \| + \sigma \|W_n \|  + O(\alpha^2) 
\end{eqnarray}
$\square$


\bibliographystyle{apalike}

\end{document}